\documentclass[12pt]{amsart}

\usepackage{amssymb}
\usepackage{amsthm}
\usepackage{amsmath}
\usepackage{graphicx}
\usepackage[only,llbracket,rrbracket]{stmaryrd}
\SetSymbolFont{stmry}{bold}{U}{stmry}{m}{n}
\usepackage{mathrsfs}
\usepackage{caption}
\usepackage[usenames,dvipsnames]{xcolor}
\usepackage{bm}
\usepackage{xspace}
\usepackage{enumitem}
\setlist[description]{leftmargin=\parindent, labelindent=\parindent}
\usepackage[normalem]{ulem}
\usepackage{outlines}
\usepackage{marginnote}
\usepackage{hyperref}
\usepackage{scalerel}
\usepackage{geometry}
\DeclareMathAlphabet{\mathpzc}{OT1}{pzc}{m}{it}

\newcommand{\map}[1]{\xrightarrow{#1}}

\newcommand{\bF}{\mathbb F}

\newcommand{\bz}{\mathbb Z}

\newcommand{\br}{\mathbb R}

\newcommand{\al}{\alpha}
\newcommand{\be}{\beta}
\newcommand{\ga}{\gamma}

\newcommand{\boldA}{\mathbf{A}}
\newcommand{\boldB}{\mathbf{B}}

\newcommand{\Ozsvath}{Ozsv{\'a}th\xspace}
\newcommand{\Szabo}{Szab{\'o}\xspace}

\DeclareMathOperator{\Sym}{Sym}

\DeclareMathOperator{\Hom}{Hom}
\DeclareMathOperator{\im}{Im}

\DeclareMathOperator{\Id}{Id}

\DeclareMathOperator{\tildecob}{\scalerel*{\widetilde{\mat}(\cob^3)}{(\cob^3)}}

\DeclareMathOperator{\CKh}{CKh}
\DeclareMathOperator{\Kh}{Kh}
\DeclareMathOperator{\rCKh}{\widetilde{CKh}}
\DeclareMathOperator{\rKh}{\widetilde{Kh}}
\DeclareMathOperator{\Sz}{Sz}
\DeclareMathOperator{\CSz}{CSz}

\DeclareMathOperator{\Lee}{Kh_{Lee}}
\DeclareMathOperator{\CBn}{\CKh_{BN}}
\DeclareMathOperator{\Bn}{\Kh_{BN}}

\newcommand{\HF}{\widehat{HF}}
\newcommand{\CF}{\widehat{CF}}

\DeclareMathOperator{\diagram}{\mathcal{D}}

\def\co{\colon\thinspace}

\newcommand{\cmap}{\mathscr{F}}
\newcommand{\decos}{\mathfrak{t}}
\newcommand{\Decos}{\mathbf{T}}
\newcommand{\config}{\mathscr{C}}

\newcommand{\kft}{\mathcal{A}}
\newcommand{\hkft}{\mathcal{K}}

\newcommand{\BN}[1]{\llbracket #1 \rrbracket}
\DeclareMathOperator{\cob}{\mathpzc{Cob}}
\DeclareMathOperator{\mat}{\mathpzc{Mat}}
\DeclareMathOperator{\matcob}{\mathpzc{Mat}(\cob^3)}

\DeclareMathOperator{\diag}{\mathpzc{Diag}}
\DeclareMathOperator{\Link}{\mathpzc{Link}}
\DeclareMathOperator{\Planar}{\mathpzc{Planar}}
\DeclareMathOperator{\Kom}{\mathpzc{Kom}}

\newcommand{\branched}{\mathpzc{Br}}

\DeclareMathOperator{\cone}{cone}
\newcommand{\handle}[1][]{\mathfrak{h}_{ #1 }}

\newtheorem{Thm}{Theorem}[section]
\newtheorem*{Thm*}{Theorem}
\newtheorem{Prop}[Thm]{Proposition}
\newtheorem*{Prop*}{Proposition}
\newtheorem{Lem}[Thm]{Lemma}
\newtheorem*{Lem*}{Lemma}
\newtheorem{Cor}[Thm]{Corollary}
\newtheorem*{Conj}{Conjecture}
\newtheorem{innercustomthm}{Theorem}

 \newcommand{\thistheoremname}{}
\newtheorem*{genericthm*}{\thistheoremname}
\newenvironment{namedthm}[1]
  {\renewcommand{\thistheoremname}{#1}%
   \begin{genericthm*}}
  {\end{genericthm*}}

\theoremstyle{definition}
\newtheorem{Def}[Thm]{Definition}
\newtheorem*{Def*}{Definition}

\theoremstyle{remark}
\newtheorem{Rem}[Thm]{Remark}
\newtheorem{Rem*}{Remark}
\newcounter{Eg}

\newtheorem*{Eg*}{Example}

\begin{document}
\title[More functoriality]{Strong Khovanov-Floer theories and functoriality}
\author{Adam Saltz}
\address{University of Georgia}
\email{adam.saltz@uga.edu}
\date{\today}
\thanks{This material is based upon work supported by the National Science Foundation under Grant No. 1664567.}

\begin{abstract}
		We provide a unified framework for proving Reidemeister-invariance and functoriality for a wide range of link homology theories.  These include Lee homology, Heegaard Floer homology of branched double covers, singular instanton homology, and \Szabo's geometric link homology theory.  We follow Baldwin, Hedden, and Lobb (arXiv:1509.04691) in leveraging the relationships between these theories and Khovanov homology.  We obtain stronger functoriality results by avoiding spectral sequences and instead showing that each theory factors through Bar-Natan's cobordism-theoretic link homology theory.  
\end{abstract}
\maketitle

\section{Introduction}\label{sec:intro}~
	Over the last twenty years,~low-dimensional topologists have been treated to a feast of new homology theories.  Perhaps the most tractable is Khovanov homology, which assigns a bigraded chain complex to a diagram of a link in $S^3$ \cite{Khovanov2000}.  The homology of this chain complex, $\Kh(L)$, is a link invariant.  Lee  perturbed Khovanov's construction to obtain a theory $\Lee$ whose rank is boring as a link invariant but whose gradings were used by Rasmussen to prove Milnor's conjecture on the four-ball genus of torus knots \cite{Lee20052,Rasmussen2010}.  Computers can compute these homologies for knots with dozens of crossings, see \cite{BarNatan2007}.

	Just as ordinary homology is a functor from the category of topological spaces to the category of abelian groups, Khovanov homology can be understood as a functor from a certain category of links to the category of abelian groups.  In particular, Khovanov homology assigns homomorphisms to link cobordisms.  (A link cobordism is a smoothly embedded surface in $S^3 \times I$ with $0$-boundary $L$ and $1$-boundary $L'$.)  This property of Khovanov homology is called \emph{functoriality.}

	Khovanov recognized a recipe for constructing maps from a diagrammatic presentation of a link cobordism but was initially unable to show that the map did not depend (up to homotopy) on the presentation.  Jacobsson \cite{Jacobsson2004} proved it using carefully organized computations, and then Khovanov \cite{Khovanov2006} gave a more algebraic proof.  Bar-Natan \cite{BarNatan2005} reproved the theorem by recasting the entire construction in terms of cobordisms.  Because of this, his proof applies to Lee's theory and other variants.

	From a totally different direction come the gauge-theoretic and symplectic invariants of three- and four-manifolds.  These include monopole Floer homology, Heegaard Floer homology, plane Floer homology, and instanton Floer homology~\cite{KronheimerMrowka,Ozsvath2004,Daemi2015,KronheimerMrowka2011}.  Each is an invariant of three-manifolds with myriad applications.  They too are functorial: a four-manifold with boundary $-Y \coprod Y'$ induces a map from the Floer homologies of $Y$ to those of $Y'$.

	The double cover of $S^3$ branched along $L$ is a three-manifold whose properties mirror those of $L$; see Montesinos' classical work \cite{Montesinos}.  The Floer homologies of $\Sigma(L)$ are link invariants.  A link cobordism from $L$ to $L'$ lifts to a four-dimensional cobordism from $\Sigma(L)$ to $\Sigma(L')$.  But it is not clear \emph{a priori} that (say) $\HF(\Sigma(L))$ is functorial as a link invariant: one needs to show that the four-dimensional cobordism maps respect relations between link cobordisms.

	The last homology theory considered here is \Szabo's ``geometric link homology,'' \cite{Szabo2015}, which we will just call ``\Szabo homology.''  We write $\CSz(\diagram)$ for the \Szabo chain complex of a link diagram and $\Sz(L)$ for \Szabo homology.  This theory is combinatorial in the same way as Khovanov and Lee's theories are.  On the other hand, Seed \cite{Seed2011} produced extensive numerical evidence that $\Sz(L) \simeq \HF(\Sigma(-L))$.  So \Szabo's theory conjecturally has a foot in both worlds.

	These invariants are connected by the following Arch-Theorem.

        \begin{namedthm}{Arch-Theorem}
		Write $-L$ for the mirror of $L$.  Let $H(L)$ be any of the following groups:
		\begin{itemize}
			\item the Heegaard Floer homology of $\Sigma(-L)$ \cite{Ozsvath2005b}
			\item the monopole Floer homology of $\Sigma(-L)$ \cite{Bloom2011}
			\item the plane Floer homology of $\Sigma(-L)$ \cite{Daemi2015}
			\item the singular instanton homology of $L$ \cite{Kronheimer2011}
			\item the framed instanton homology of $L$ \cite{Scaduto2015}
			\item the Lee-Bar-Natan\footnote{Lee and Bar-Natan each perturbations of Khovanov homology.  They are twist-equivalent \cite{Khovanov20062} over any ring in which $2$ is invertible, but they are not equivalent over $\bF$.  We will work with Bar-Natan's theory.  To avoid confusion between this homology theory and Bar-Natan's cobordism-theoretic construction, we refer to the former as \emph{Lee-Bar-Natan homology}.} homology of $L$ \cite{Rasmussen2010}
			\item the \Szabo homology of $L$ \cite{Szabo2015}
		\end{itemize}
		For each diagram of $L$ there is a spectral sequence which converges to $H(L)$ and whose second page is isomorphic to $\Kh(L)$.  (In some instances one should replace $\Kh(L)$ with a reduced version.)
		\end{namedthm}

	Beyond their inherent interest, these spectral sequences have many topological applications, see the introduction of \cite{BaldwinHeddenLobb2015}.  For example, the only extant proof that Khovanov homology detects the unknot uses the spectral sequence to singular instanton.  What is the common structure behind all of this?  

	Baldwin, Hedden, and Lobb \cite{BaldwinHeddenLobb2015} define a \emph{Khovanov-Floer theory} to be a gadget $\kft$ which assigns to a link diagram a spectral sequence whose second page is isomorphic to Khovanov homology.  The theory must come with some additional data.  For example,  if $\diagram'$ is the result of a one-handle attachment on $\diagram$, there is a map of spectral sequences $F^i \co \kft(\diagram) \to \kft(\diagram')$ so that $F^2$ agrees with the one-handle attachment map on Khovanov homology.  Using this barebones structure, they prove the following.

    \begin{Thm*}[\cite{BaldwinHeddenLobb2015}]
		Let $\kft$ be a Khovanov-Floer theory.
		\begin{enumerate}
			\item $\kft$ is a link invariant: if $\diagram'$ is obtained from $\diagram$ by a sequence of Reidemeister moves, then there is an isomorphism\footnote{An isomorphism of spectral sequences is a map of spectral sequences which is an isomorphism after the first page.} of spectral sequences $\kft(\diagram) \cong \kft(\diagram')$.  
			\item $\kft$ is functorial: a movie presentation $M(\Sigma)$ of a link cobordism $\Sigma$ from $L$ to $L'$ induces a map of spectral sequences
			\[
				F^i_{M(\Sigma)} \co \kft(\diagram) \to \kft(\diagram').
			\]
			If $M'(\Sigma)$ is some other movie presentation of $\Sigma$, then $F^i_{M(\Sigma)} = F^i_{M'(\Sigma)}$ for $i \geq 2$.		
			\item Every theory in the Arch-Theorem is a Khovanov-Floer theory.
		\end{enumerate}
	\end{Thm*}
	Ideally, one could use this definition to classify or axiomatize ``homologies which admit spectral sequences from Khovanov homology'' in the spirit of Eilenberg-Steenrod for ordinary homology and Khovanov for link homologies constructed from rank two Frobenius algebras \cite{Khovanov20062}.

	Unfortunately, the construction is limited by its reliance on spectral sequence techniques.  To see this, consider a movie presentation $M(\Sigma)$ of a link cobordism $\Sigma$.  It induces a filtered map of \Szabo chain complexes
    \[
		F_{M(\Sigma)} \co \CSz(\diagram) \to \CSz(\diagram').
	\] 
    ``\Szabo homology is functorial'' should mean that the induced map
    \[
		(F_{M(\Sigma)})_* \co \Sz(\diagram) \to \Sz(\diagram')
	\] 
    does not depend on the presentation of $\Sigma$.  In the Khovanov-Floer setting, $M(\Sigma)$ induces a map of spectral sequences $F^i_{M(\Sigma)}$.  This map converges\footnote{We only consider bounded complexes with coefficients in a field, so there's no question of convergence.} to a different map, $F^\infty_{\Sigma}$, and this is the map studied in \cite{BaldwinHeddenLobb2015}.  These maps may be quite different, even in simple theories; see the end of Section \ref{subsec:LeeTurner}.  So Baldwin, Hedden, and Lobb show that the \emph{Khovanov-Floer theory associated to \Szabo homology} is functorial, but not necessarily \Szabo homology itself.

	Our goal in this paper is to formulate a notion of Khovanov-Floer theory which does not rely on spectral sequences.  Instead, a strong Khovanov-Floer theory assigns a filtered chain complex to a link.  To a zero-, one-, or two-handle attachment the theory assigns a chain map.  These maps are required to satisfy some relations which are fundamentally three-dimensional: they concern cobordisms of planar link diagrams.  If the theory satisfies the technical condition of being \emph{conic}, then its filtered chain homotopy type is a functorial link invariant link invariant in the usual sense.  This is our main result.

	\begin{innercustomthm}
		Let $\hkft$ be a strong Khovanov-Floer theory.
		\begin{enumerate}
			\item $\hkft$ is a link invariant: if $\diagram'$ is obtained from $\diagram$ by a sequence of Reidemeister moves, then $\hkft(\diagram) \simeq \hkft(\diagram')$.
			\item $\hkft$ is functorial: a movie presentation $M(\Sigma)$ of a link cobordism $\Sigma$ from $L$ to $L'$ induces a chain map
			\[
				F_{M(\Sigma)} \co \hkft(\diagram) \to \hkft(\diagram').
			\]
			If $M'(\Sigma)$ is some other movie representation of $(\Sigma)$, then $F_{M_\Sigma} \simeq F_{M'_\Sigma}$.		
			\item Lee-Bar-Natan homology, \Szabo homology, Heegaard Floer homology, and singular instanton homology yield strong Khovanov-Floer theories.
		\end{enumerate}
	\end{innercustomthm}

	\begin{innercustomthm}
		\Szabo homology, Heegaard Floer homology of branched double covers, and singular instanton link homology are functorial link invariants.\footnote{Functoriality of Lee-Bar-Natan homology has been known for some time.}
	\end{innercustomthm}

	We expect that every theory in the Arch-Theorem yields a strong Khovanov-Floer theory.

	Our central technical result is Proposition \ref{prop:mainTechnical}, which states roughly that every conic, strong Khovanov-Floer theory factors through Bar-Natan's cobordism-theoretic construction of link homology \cite{BarNatan2005}.  Khovanov homology is usually defined by applying some linear algebraic machinery to the cube of resolutions of a link diagram.  Bar-Natan instead formulates a category in which the resolutions are objects and the saddle cobordisms between them are morphisms.  The cube of resolutions is a complex $\BN{\diagram}$ in this category.  The chain homotopy type of this complex is a link invariant.

	Let $\diagram$ and $\diagram'$ be link diagrams for a link $L$.  Write $\BN{\diagram}$ and $\BN{\diagram'}$ for their cube of resolutions complex.  Suppose that $f \co \BN{\diagram} \to \BN{\diagram'}$ is one direction of a chain homotopy equivalence.  This means that $f$ is a `map' comprising a formal sum of cobordisms.  The observation which set off this work was that these maps, properly interpreted, also furnish a chain homotopy equivalence for \Szabo homology.

	One cannot discuss functoriality without paying some attention to naturality.  Each link homology theory (aside from Bar-Natan-Lee homology and Khovanov homology) requires some auxiliary input in addition to a link diagram.  This is why Theorem 2 is limited to only a few examples.  We expect that every theory in the Arch-Theorem is a conic, strong Khovanov-Floer theory, and that the proof for each should be very similar to the proof for singular instanton homology, but we have not checked the necessary naturality conditions.

	\subsection{Future directions}
		Strong Khovanov-Floer theories should be the right framework to extend and universalize other results in Khovanov homology.  For example, many Khovanov-Floer theories include invariants of contact structures and transverse links.  Baldwin and Plamenevskaya \cite{Baldwin2010}, following Roberts \cite{Roberts2013}, showed that these invariants in Khovanov homology and Heegaard Floer homology are connected.  In Bar-Natan's framework we can define a `universal' transverse invariant which we believe will be useful in understanding the connections between all of these invariants.

		Stephan Wehrli showed that, over $\bz/2\bz$, Khovanov homology is invariant under knot mutation.  His proof runs through Bar-Natan's framework, but uses more structure than we construct here.  Nevertheless, we believe it can be adapted to prove the following.
		\begin{Conj}
			Every strong Khovanov-Floer theory over $\bz/2\bz$ is invariant under knot mutation.
		\end{Conj}
		Lambert-Cole \cite{2017arXiv170100880L} defined \emph{extended Khovanov-Floer theories} and showed that they are mutation invariant.  It would be interesting to adapt his proof to strong Khovanov-Floer theories. 

		Lin \cite{Lin2016} has shown that there is a spectral sequence from (reduced) Bar-Natan-Lee homology to involutive monopole Floer homology.  There should be no difficulty in adapting the definition of strong Khovanov-Floer theory to include this case.

		Lastly, John Baldwin challenged the author to formulate \Szabo link homology theory in cobordism-theoretic terms.  This paper does not meet that challenge: we cannot describe \Szabo's configuration maps purely in terms of cobordisms.  But we believe that these results suggest a path forward to better understand \Szabo homology along with other Khovanov-Floer theories.  Putyra \cite{Putyra2014} unified even and odd Khovanov homology using the theory of \emph{framed functions}, an enrichment of traditional Morse (really, Cerf) theory.  These framings are nearly identical to \Szabo's decorations.  We conjecture that \Szabo's theory could be given a purely Morse-theoretic definition using framed functions.  If \Szabo and Seed's conjecture holds, this would give a purely two-dimensional definition for Heegaard Floer homology of branched double covers.   An optimist might conjecture every conic, strong Khovanov-Floer theory can be understood via two-dimensional Morse theory by working with an Cerf-enhanced version of Bar-Natan's construction.

	\subsection{Organization and conventions}

		In Section \ref{sec:diagrams}, we discuss the connections between links, link diagrams, cobordisms, and diagrammatic cobordisms.  In Section \ref{sec:khovanovFloerTheories} we review BHL's definition and define strong Khovanov-Floer theories.    Section \ref{sec:cones} is a review of Bar-Natan's cobordism-theoretic framework.  Section \ref{sec:cobordismsI} is the technical core of the paper in which we prove the first two parts of the main theorem.  We prove that the examples in the Arch-Theorem are actually strong Khovanov-Floer theories in Sections \ref{sec:examples} and \ref{sec:hardExamples}.

		$\mathbb{F}$ stands for the field $\bz/2\bz$.  Tensor products are assumed to be over $\mathbb{F}$.  For other coefficient systems, see the remark at the end of section \ref{subsec:hkft}.

\subsection*{Acknowledgments} 
	This paper owes its existence to John Baldwin in several ways: first, for pushing me to study the functoriality of \Szabo homology while I was his student; second, for his careful demolition of arguments in \cite{Baldwin2018} which showed the necessity of the present argument; and third, for his part in authoring \cite{BaldwinHeddenLobb2015}.

	Thanks to David Gay and Gordana Mati{\`c} for many inspirational discussion.  Andrew Lobb and David Rose made a valuable corrections to an earlier version.  I am also thankful to Matt Hedden, Peter Lambert-Cole, Tye Lidman, Mike Usher, Mike Willis for helpful conversations.

\section{Diagrams and cobordisms}\label{sec:diagrams}

	\subsection{Diagrammatic link cobordisms and embedded cobordisms}
	
		One of the first steps towards showing Khovanov homology is a link invariant is to show that it is invariant under planar isotopy.  One constructs a chain homotopy equivalence between the  Khovanov chain complexes of two diagrams which differ by a planar isotopy.  These homotopy equivalences are not natural: there are planar isotopies whose initial and final diagrams are identical but which are not assigned the identity map.  This means that any functoriality statement for Khovanov homology -- or any flavor of Khovanov-Floer theory -- must make careful distinctions between diagrams, links, and cobordisms of each.

		\begin{Def}
			Let $\Link$ be the category of links in $\br^3$.  Its objects are smooth links in $\br^3$ and its morphisms are (proper) isotopy classes of collared, smoothly embedded link cobordisms in $\br^3 \times I$.  These isotopies are always taken to fix a neighborhood of the boundary pointwise.  
		\end{Def}

		Carter and Saito \cite{CarterSaito1993} set up a combinatorial model for $\Link$.

\begin{figure}
			\begin{center}
			\includegraphics[width=.75\linewidth]{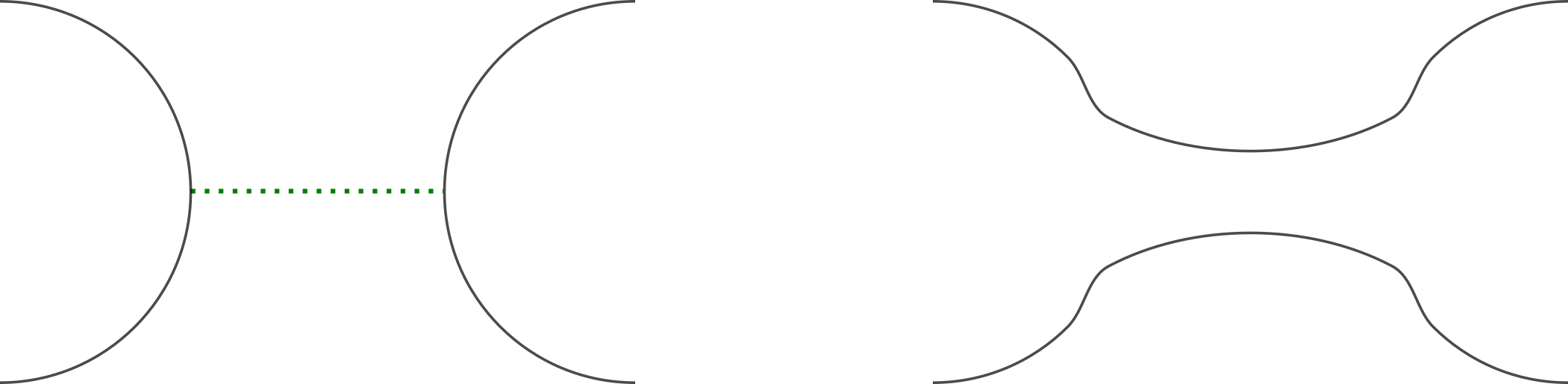}
			\caption{Before and after a diagrammatic one-handle attachment along the green, dotted arc.}
			\label{fig:diagrammatic1Handle}
			\end{center}
		\end{figure}

		\begin{Def}\label{def:diagHandle}
			Let $\diagram$ be a link diagram.  A \emph{diagrammatic handle attachment} on $\diagram$ is one of the three following operations:
			\begin{description}[labelwidth=3cm]
				\item[Zero-handle attachment] add a disjoint circle to $\diagram$.
				\item[Two-handle attachment] remove a crossingless, disjoint circle from $\diagram$.
				\item[One-handle attachment] Let $\gamma$ be an arc embedded in $\br^2$ so that $\partial \gamma \subset \diagram$ and $\gamma$ is otherwise disjoint from $\diagram$.  A one-handle attachment along $\gamma$ is the local operation represented in Figure \ref{fig:diagrammatic1Handle}.
			\end{description}
			The \emph{support} of a zero- or two-handle attachment is a small neighborhood of the relevant circle.  The support of a one-handle attachment along $\gamma$ is a small neighborhood of $\gamma$.
		\end{Def}		

		Let $L$ and $L'$ be links and let $\Sigma$ be a generic representative of an element in $\Hom_{\Link}(L,L')$.  Write $t$ for the last coordinate in $\br^3 \times I$.  Write $\Sigma_t = \Sigma \cap (\br^3 \times \{t\})$.  For almost all $t$, $\Sigma_t$ is a link in $\br^3$ whose projection $\diagram_t$ to the $xy$-plane is a regular diagram.  Call these values of $t$ \emph{regular}.  If the interval $[t_0,t_1]$ consists only of regular values then the cobordism $\Sigma_{t_0} \to \Sigma_{t_1}$ is a cylinder and the cobordism $\diagram_{t_0} \to \diagram_{t_1}$ can be represented as a planar isotopy.  If $t_2$ is the only singular value in $[t_0, t_1]$ then $\diagram_{t_1}$ is obtained from $\diagram_{t_0}$ by either a diagrammatic handle attachment or a Reidemeister move.  Let
		\[
			[\diagram_1, \ldots, \diagram_n]
		\]
		be a sequence of diagrams so that $\diagram_{i+1}$ is obtained from $\diagram_i$ by an \emph{elementary diagrammatic cobordism:} a planar isotopy, a single diagrammatic handle attachment, or a single Reidemeister move.  Such a sequence is called a \emph{movie}.  The constituent diagrams are called \emph{frames}.  

\begin{figure}
			\begin{centering}
				\includegraphics[width=.5\linewidth]{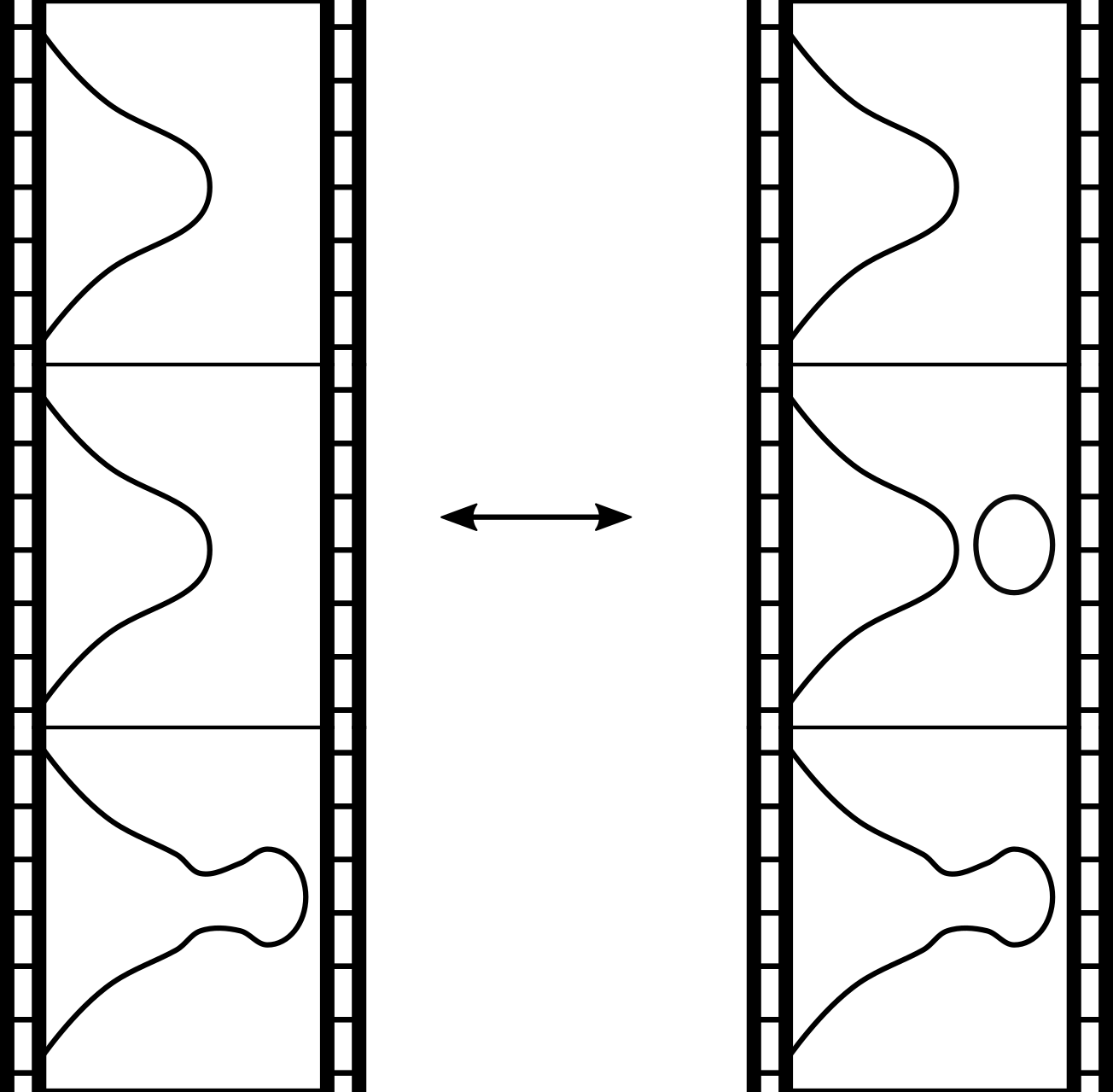}
				\caption{Movie move 15.  These are local pictures of two diagrammatic cobordisms which are assumed to be equal elsewhere.  The embedded cobordisms represented by these movies are isotopic.  (The film style on this diagram is due to Carter and Saito.)}
				\label{fig:movieMove16}
			\end{centering}
		\end{figure}

		It is not too hard to see that a movie describes a link cobordism.  In fact, any link cobordism admits a movie description.  One can alter a movie in a way which preserves the isotopy class of the presented surface.  An example is given in Figure \ref{fig:movieMove16}.  Formally, a \emph{movie move} is a rule for replacing one collection of frames with another so that the isotopy class of the presented surface does not change.  Carter and Saito identify a set of movie moves which are sufficient to connect any two movies which represent the same isotopy class of surface.

		\begin{Thm}[\cite{CarterSaito1993}]\label{ref:CarterSaito} 
			Let $L$ and $L'$ be links.  Every $\Sigma \in \Hom_{\Link}(L,L')$ has a representative which can be described by a movie.  Any two movies for the same cobordism can be related by applying a sequence of the movie moves.
		\end{Thm}

		The upshot is that any invariant of surfaces defined by slicing them into movies must be invariant under the movie moves.

		\begin{Def}\label{def:diag}
			Let $\diag$ be the category whose objects are link diagrams and whose morphisms are movies modulo the movie moves.  That is, $\Hom_{\diag}(\diagram,\diagram')$ is the set of movies with initial frame $\diagram$ and final frame $\diagram'$ modulo the relation $M \sim M'$ if $M$ and $M'$ can be connected by a sequence of movie moves.
		\end{Def}

		\begin{Prop}[\cite{BaldwinHeddenLobb2015}]\label{prop:equivDiagCob}
			Let $\pi \co \br^3 \to \br^2$ be the projection to the $xy$-plane.  Let $\mathcal{L}_0$ be the set of links $L$ so that $\pi(L)$ is immersed.  

			There is an equivalence of categories $\Pi \co \Link \equiv \diag$ which sends $L \in \mathcal{L}_0$ to $\pi(L)$.
		\end{Prop}
		
		The equivalence is given by fixing an isotopy between each object of $\Link$ and an object which projects to a regular immersed curve.  These isotopies may be viewed as topological cylinders in $\br^3 \times I$ and therefore as morphisms in $\Link$.  This data defines a functor (and equivalence) from $\Link$ to $\diag$.

		Any two choices of isotopies induce isomorphic functors.  So one may choose a set of isotopies which fix $\mathcal{L}_0$.  (In fact, one may choose isotopies which are the identity on any link which already projects to a regular immersion.)  There are many equivalences.  Fix one and call it $\Pi$.

		The two movie moves which play a major role in this paper are movie move 15 (shown in Figure\footnote{In fact this picture encodes two moves: it can be read from top to bottom or bottom to top.  We will refer to this package of moves as ``movie move 15.''} \ref{fig:movieMove16}) and the \emph{distant handleswap}.  The latter is the operation of swapping the order of handle attachments with disjoint supports.  Both of these moves have three-dimensional interpretations: movie move 15 is the cinematic version of the addition or removal of a pair of canceling handles for a surface in $\br^2 \times I$.  Distant handleswap is an isotopy whose support has two disjoint components: the supports of the two swapped handles.

	\subsection{The categories of planar diagrams and planar circles}

		\begin{Def}
			Let $\Planar$ be the category whose objects are embedded collections of circles in $\br^2$ and whose morphisms are diagrammatic handle attachments and planar isotopies modulo movie move 15 and distant handleswaps.  $\Planar$ is the \emph{category of planar diagrams}.
		\end{Def}

		\begin{Def}
			Let $\cob^3$ be the category whose objects are embedded collections of circles in $\br^2$ and whose morphisms are cobordisms in $\mathbb{R}^2 \times I$.  Two cobordisms are considered equal if they are related by a boundary-preserving, proper isotopy.
		\end{Def}

		The following lemma is a three-dimensional version of Carter and Saito's work.  Consequently, it is much simpler.

		\begin{Lem}\label{lem:movies2}
			Let $O$ and $O'$ be objects of $\cob^3$.  Let $\Sigma \in \Hom_{\cob^3}(O,O')$.  Then $\Sigma$ admits a movie representative $M(\Sigma)$ with only crossingless diagrams.  Any two movie representatives for $\Sigma$ from this construction differ by swapping distant saddles or movie move 15.  

			Any movie $M$ with first frame $O$ and last frame $O'$ whose frames differ only by planar isotopy and diagrammatic handle attachments describes a cobordism $\Sigma_M \in \Hom_{\cob^3}(O,O')$ so that $M(\Sigma_M) = M$.  If $M$ and $M'$ are related by distant handleswaps and movie move 15, then $\Sigma_M = \Sigma_{M'}$. 
		\end{Lem}

		\begin{proof}
			The first part of the theorem may be restated as follows: let $\Sigma$ and $\Sigma'$ be isotopic surfaces in $\br^2 \times I$ so that the natural height function is Morse on both.  Assume also that both height functions are \emph{separating:} no two critical points share a critical level.  Then there is a one-parameter family of embeddings (i.e. an isotopy) $f_t$ from $\Sigma$ to $\Sigma'$ which is Morse and separating except at finitely many points.  At those points, either two critical points share a critical level or there is a birth-death singularity.  Diagrammatically, these are handleswaps and applications of movie move 15.

			Let $\tilde{g}_t$ be an isotopy from $\Sigma$ to $\Sigma'$.  It is standard Cerf theory that there is a function $h \co \Sigma \to [0,1]$ so that, for any $\epsilon > 0$, the function $g_t = \tilde{g}_t + \epsilon h$ is separating and Morse except for finitely many handleswaps and birth-death singularities.  If $\epsilon$ is small enough, then $g_t$ is an embedding for all $t$.  Therefore $g_t$ is an isotopy whose cinematic presentation is by handleswaps and movie move 15.

			A movie describes a surface by thickening frames into cylinders and attaching handles across diagrammatic handle attachments.  The last statement is clear.
		\end{proof}

		Let $\Pi \co \Link \equiv \diag$ be the equivalence from the discussion following Proposition \ref{prop:equivDiagCob}.  Let $P$ be the plane of projection for $\Pi$.  Choose a level-preserving embedding of $\br^2 \times I$ into $\br^3 \times I$ so that $\br^2 \times \{0\}$ is mapped to $P \times \{0\}$ and $\br^2 \times \{1\}$ is mapped to $P \times \{1\}$.  This induces an inclusion functor $\iota \co \cob^3 \to \Link$.

		\begin{Prop}\label{prop:equivDiag2Cob2}
			The equivalence $\Pi \co \Link \equiv \diag$ restricts to an equivalence $\Pi \circ \iota \co \cob^3 \equiv \Planar$.
		\end{Prop}

		\begin{proof}
			Let $O$ be an object of $\cob^3$.  It's clear that $(\Pi \circ \iota)(\gamma)$ is an object of $\Planar$.   It's also obvious that $(\Pi \circ \iota)$ is essentially surjective.  The fullness and faithfulness of $\Pi \circ \iota$ follows from Lemma \ref{lem:movies2}.
		\end{proof}

		Proposition \ref{prop:equivDiag2Cob2} justifies our future conflation of cobordisms and diagrammatic cobordisms.
\section{Strong Khovanov-Floer theories}\label{sec:khovanovFloerTheories}
	\subsection{Khovanov-Floer theories}\label{subsec:kfts}
		We recall Baldwin, Hedden, and Lobb's definition of Khovanov-Floer theory.  To get a reasonable category of ``filtered complexes whose spectral sequence has second page isomorphic to $\Kh(L)$'' one has to keep track of the isomorphisms.

		\begin{Def}
			Let $V$ be a graded vector space.  A \emph{$V$-complex} is a pair $(C,q)$ where $C$ is a filtered chain complex and $q$ is an isomorphism from $V$ to $E_2(C)$. A map of $V$-complexes is a filtered chain map.  A map $f \co (C,q) \to (C',q')$ induces a map $E_2(f) \co E_2(C) \to E_2(C')$ which may be written as $q' \circ T \circ q^{-1}$ for some map $T \co V \to W$.  Say that \emph{$f$ agrees with $T$ on $E_2$.}  If $V = W$ and $f$ agrees with the identity map then call $f$ a \emph{quasi-isomorphism}.  
		\end{Def} 

		\begin{Def}\label{def:kft}
			A \emph{Khovanov-Floer theory} $\kft$ is a rule which assigns to every link diagram $\diagram$ a quasi-isomorphism class of $\Kh(\diagram)$-complexes $\kft(\diagram)$ which satisfies the following rules.
			\begin{enumerate}
				\item Suppose that $\diagram'$ is obtained from $\diagram$ by a planar isotopy.  Then there is a morphism
				\[
					\kft(\diagram) \to \kft(\diagram')
				\]
				which agrees on $E_2$ with the induced map from $\Kh(\diagram)$ to $\Kh(\diagram')$.
				\item Suppose that $\diagram'$ is obtained from $\diagram$ by a diagrammatic $1$-handle attachment.  Then there is a morphism
				\[
					\kft(\diagram) \to \kft(\diagram')
				\]
				which agrees on $E_2$ with the induced map from $\Kh(\diagram)$ to $\Kh(\diagram')$.
				\item Write $\diagram \cup \diagram'$ for the disjoint union of the diagrams $\diagram$ and $\diagram'$.  There is a morphism
				\[
					\kft(\diagram \cup \diagram') \to \kft(\diagram) \otimes \kft(\diagram') 
				\]
				which agrees on $E_2$ with the standard isomorphism
				\[
					\Kh(\diagram \cup \diagram') \to \Kh(\diagram) \otimes \Kh(\diagram').
				\]
				\item If $\diagram$ is a diagram of the unlink, then $E_2(\kft(\diagram)) \cong E_{\infty}(\kft(\diagram))$.
			\end{enumerate}
		\end{Def}
	
		\begin{Def}\label{def:weakFunctoriality}
			A Khovanov-Floer theory is \emph{functorial} if, for a movie $M \co \diagram \to \diagram'$, there is a filtered chain map
			\[
				F_M \co \kft(\diagram) \to \kft(\diagram)
			\]
			which agrees which agrees on $E_2$ with the usual map on Khovanov homology \cite{Jacobsson2004}.  Equivalently, link cobordisms induce maps of spectral sequences which are well-defined (i.e. do not depend on the presentation) after the first page.
		\end{Def}

		\begin{Thm*}[Baldwin, Hedden, Lobb \cite{BaldwinHeddenLobb2015}] 
			Every Khovanov-Floer theory is functorial.
		\end{Thm*}

		As noted in the introductions of \cite{BaldwinHeddenLobb2015} and the present work, this functoriality is not the same as functoriality of a link homology theory.

	\subsection{Transitive systems of chain complexes}

		This section establishes a naturality framework for working with chain homotopy classes of maps.  The first section of Juhasz and Thurston \cite{Juhasz2012} is a good introduction to naturality issues in homology theories.

		\begin{Def}
			A \emph{transitive system of filtered chain complexes up to homotopy} (or \emph{transitive system} for short) is a collection of filtered chain complexes $\{C_\al\}$ and a collection of filtered chain maps $\psi_{\al}^{\be} \co C_{\al} \to C_{\be}$ for all $\al$ and $\be$ in some index set $I$ so that
			\begin{itemize}\setlength{\itemsep}{0em}
				\item $\psi_\al^\al \simeq \Id_{\al}$ where $\Id_{\al}$ is the identity map on $C_\al$.
				\item $\psi_\be^\ga \circ \psi_\al^\be \simeq \psi_{\al}^{\ga}$.
			\end{itemize}

			A transitive system can be viewed as a diagram in the homotopy category of chain complexes.  The \emph{canonical representative} of a transitive system is its inverse limit.
		\end{Def}

		The following lemma follows immediately from the fact that the inverse limit is a functor.

		\begin{Lem}\label{lem:transitiveMaps}
			Let $C = \{C_\al\}$ and $C' = \{C'_\al\}$ be transitive systems.  Let $f$ be a map of transitive systems, i.e. a collection of filtered chain maps $f_\al^\ga \co C_\al \to C'_\ga$ so that
			\[
				f_\be^\ga \circ \psi_{\al}^{\be} \simeq \psi_{\delta}^\ga \circ f_{\al}^{\delta}
			\]
			holds for any indices $\al$, $\be$, $\ga$, and $\delta$.  Then $f$ defines a filtered map between the canonical representatives of $C$ and $C'$.  If the maps $f_{\be}^\ga$ are chain homotopy equivalences then $f$ is a chain homotopy equivalence.
		\end{Lem}

		\begin{Lem}\label{lem:naturality}
			Let $C$ and $C'$ be transitive systems with maps $\psi$ and $\psi'$.  Let $I_{C}$ and $I_{C'}$ be the index sets for $C$ and $C'$.  Suppose there that there is a surjective function $\phi \co I_{C} \to I_{C'}$ and a collection of maps
			\[
				f = \{f_\al^{\phi(\al)} | \al \in I_{C}\}
			\]
			so that 
			\[
				\psi'^{\phi(\be)}_{\phi(\al)} \circ f_{\al}^{\phi(\al)} \simeq f_{\be}^{\phi(\be)} \circ \psi_{\al}^{\be}.
			\]  
			Then $f$ extends to a map of transitive systems and therefore of canonical representatives.  The extension is unique up to homotopy.
		\end{Lem}
		\begin{proof}
			Let $\be \in I_{C'}$.  Let $\gamma \in \phi^{-1}(\be)$.  Note that $f_{\ga}^{\be}$ is determined by hypothesis.  For $\al \notin \phi^{-1}(\be)$, define
			\[
				f_{\al}^{\be} = f_{\ga}^{\be} \circ \psi_{\al}^{\ga}.
			\]
			Suppose that $\gamma'$ is some other element of $\phi^{-1}(\be)$.  Then 
			\[
				f_{\al}^{\be} = f_{\ga}^{\be} \circ \psi_{\al}^{\ga} \simeq f_{\ga'}^{\be} \circ \psi^{\ga'}_{\ga} \circ \psi^\ga_{\ga'} \psi_{\al}^{\ga'} = f_{\ga'}^{\be} \circ \psi_{\al}^{\ga'}.
			\]
			Therefore the definition does not depend on the choice of an element of $\phi^{-1}(\be)$.

			Next we check that this actually defines a map of systems.  Suppose that $\al$ is in the complement of $\phi^{-1}(\be)$ in $I_C$.  Let $\delta \in I_{C}$.  In the next set of equations we write $\phi^{-1}(\be)$ to stand for a fixed element of that set.  Then
			\begin{align*}
				f_{\al}^{\be} \circ \psi_{\delta}^{\al} &\simeq (f_{\phi^{-1}(\be)}^{\be} \circ \psi_{\al}^{\phi^{-1}(\be)}) \circ \psi_{\delta}^{\al}  \\
				&\simeq (\psi'^{\be}_{\phi(\al)} \circ f_{\al}^{\phi(\al)}) \circ \psi^{\al}_{\delta} \\
				&=  \psi'^{\be}_{\phi(\al)} \circ (f_{\al}^{\phi(\al)} \circ \psi^{\al}_{\delta}) \\
				&\simeq \psi'^{\be}_{\phi(\al)} \circ (\psi'^{\phi(\al)}_{\phi(\delta)} \circ f^{\phi(\delta)}_{\delta}) \\
				&= (\psi'^{\be}_{\phi(\al)} \circ \psi'^{\phi(\al)}_{\phi(\delta)}) \circ f^{\phi(\delta)}_{\delta} \\
				&\simeq \psi'^\be_{\phi(\delta)} \circ f^{\phi(\delta)}_{\delta}
			\end{align*}
		as required.
		\end{proof}

		Thanks to transitivity, $\phi$ need not be surjective.

		\begin{Lem}\label{lem:naturality2}
			Suppose that either $I_C$ and $I_{C'}$ are both empty or that $I_C$ is not empty.  Then Lemma \ref{lem:naturality} holds even if $\phi$ is not surjective.
		\end{Lem}
		\begin{proof}
			Let $\be \in I_{C'} \setminus \im(\phi)$ and $\ga \in \im(\phi)$.  Let $\al' \in \phi^{-1}(\ga)$.  For any $\al \in I_{C}$, define
			\[
				f_\al^\be = \psi'^\be_{\ga} \circ f_{\al'}^{\ga} \circ \psi_{\al}^{\al'}.
			\]
			One can show that this does not depend on $\gamma$ as in the last proof.  To show that $\{f_\al^\be\}$, repeat the second part of the last proof on both sides of $\psi'^\be_{\ga} \circ f^{\ga}_{\al'} \circ \psi^{\al'}_{\al}$.
		\end{proof}

	\subsection{Strong Khovanov-Floer theories}\label{subsec:hkft}

		\begin{Def}\label{def:gkft}
			A \emph{strong Khovanov-Floer theory} $\hkft$ is a rule which assigns to a link diagram $\diagram$ and a collection of auxiliary data $A$ a filtered chain complex $\hkft(\diagram,A)$ satisfying the following conditions.
			\begin{enumerate}\setlength\itemsep{0em}
				\item\label{kft:coherence} Let $A_\al$ and $A_\be$ be two collections of auxiliary data.  Then there is a chain homotopy equivalence
				\[
					a_{\al}^\be \co \hkft(\diagram,A_\al) \to \hkft(\diagram, A_\be).
				\] 
				The collection $\{\hkft(\diagram,A_\al\}$ forms a transitive system with maps $\{a_\al^\be\}$.  Write $\hkft(\diagram)$ for the canonical representative.
				\item\label{kft:crossingless} If $\diagram$ is a crossingless diagram of the unknot, then $H(\hkft(\diagram)) \cong \Kh(\diagram)$.
	            \item\label{kft:kunneth} Let $\diagram \cup \diagram'$ be a disjoint union of diagrams.  Then 
	            \[
	            	\hkft(\diagram \cup \diagram') \simeq \hkft(\diagram) \otimes \hkft(\diagram').
	            \]
	        \end{enumerate}
	        A strong Khovanov-Floer theory also assigns maps to diagrammatic cobordisms with auxiliary data.  These maps must satisfy the following properties:
	        \begin{enumerate}[resume]
	            \item\label{kft:handles} Suppose that $\diagram'$ is obtained from $\diagram$ by a diagrammatic handle attachment.  There is a function $\phi$ from auxiliary data for $\diagram$ to auxiliary data for $\diagram'$ and a map
	            \[
	            	\handle[(A_{\al},\phi(A_{\al}),B)] \co \hkft(\diagram,A_\al) \to \hkft(\diagram',\phi(A'_{\al'})).
	            \]
	            where $B$ is some additional auxiliary data.  For fixed $B$, these maps satisfy the conditions of Lemma \ref{lem:naturality2} and therefore form a map
	            \[
	            	\handle[B] \co \hkft(\diagram) \to \hkft(\diagram').
	            \]
	            For any two sets of additional data $B$ and $B'$, $\handle[B] \simeq \handle[B']$.
	            \item\label{kft:frobenius} Let $U$ be a crossingless diagram of an unknot.  Then $\hkft(U)$ is a Frobenius algebra with operations given by the handle attachment maps.  (We will usually assume that the Frobenius algebra is $\bF[X]/(X^2)$.)
				\item\label{kft:isotopy} If $\diagram'$ is obtained from $\diagram$ by a planar isotopy, then $\hkft(\diagram)$ is homotopy equivalent to $\hkft(\diagram')$.
				\item\label{kft:handleKunneth} Let $\diagram$ be the disjoint union of two diagrams $\diagram_0$ and $\diagram_1$.  Let $\Sigma$ be a diagrammatic cobordism to $\diagram'$.  Suppose that $\diagram'$ is the disjoint union of $\diagram'_0$ and $\diagram'_1$.  Suppose that $\Sigma$ is the disjoint union of $\Sigma_0$ and $\Sigma_1$, where $\Sigma_i$ is a cobordism from $\Sigma_i$ to $\Sigma'_i$.  Then $\hkft(\Sigma) \simeq \hkft(\Sigma_0) \otimes \hkft(\Sigma_1)$.
				\item\label{kft:conditions} The handle attachment maps satisfy handleswap invariance rule and the movie move 15 relation up to homotopy.
			\end{enumerate}
		\end{Def}

		Again: $\hkft(\diagram)$ refers to a single chain complex: the canonical representative of the transitive system defined by $\diagram$.

	    The naturality issues in conditions \ref{kft:coherence} and \ref{kft:handles} are significant.  For example, to construct a Heegaard Floer complex for a three-manifold $Y$ one needs a Heegaard diagram and some analytic data.  Naturality of the complex under variations of analytic data was established by \Ozsvath and \Szabo in their foundational work on the theory.  Naturality with respect to Heegaard diagrams was proven later by Juhasz and Thurston and requires careful study of the space of all (pointed) Heegaard diagrams for a three-manifold \cite{Juhasz2012}.

	    In practice we only consider strong Khovanov-Floer theories which satisfy the following property.  Let $c$ be a crossing of a link diagram $\diagram$.  Let $\diagram_0$ be the diagram given by zero-resolving $c$.  Let $\gamma_c$ be a planar arc between the two sides of the resolution.  (Compare Figures \ref{fig:diagrammatic1Handle} and \ref{fig:resolutions}.)

	    \begin{Def}\label{def:conic}
	    	Let $\hkft$ be a strong Khovanov-Floer theory.  We say that $\hkft$ is \emph{conic} if it satisfies the following property: for any link diagram $\diagram$ and any crossing $c$ of $\diagram$,
	    	\[
					\hkft(\diagram) \simeq \cone(\handle[\gamma_c] \co \hkft(\diagram_0) \to \hkft(\diagram_1)).
			\]  
			where $\diagram_0$ and $\diagram_1$ are the $0$- and $1$-resolutions of $\diagram$ at $c$ and $\handle[\ga_c]$ is the handle attachment map at $c$.
		\end{Def}

		We will see in Section \ref{subsec:cones} that being conic is equivalent to having a cube of resolutions description.  Every strong Khovanov-Floer theory we know is conic.  Conicity and mapping cone formulas have played major roles in the application of Floer theory to low-dimensional topology, so it is not surprising to see it appear here.  It may be possible to prove that every strong Khovanov-Floer theory is conic using techniques similar to those used in \cite{BaldwinHeddenLobb2015} to show that every Khovanov-Floer theory assigns maps to zero- and two-handles.  In any case, there may be redundancy among the conditions of Definition \ref{def:gkft} and conicity.  To start, Corollary \ref{cor:conicMeansCondition} at the beginning of Section \ref{sec:examples} shows that conicity implies Condition \ref{kft:conditions}.

		\begin{Rem}\label{rem:tqft}
			This paper aims to show that Definition \ref{def:gkft} is the right way to formulate ``link homology theory with a spectral sequence from Khovanov homology.''  In light of Condition \ref{kft:frobenius} of Definition \ref{def:gkft} and the well-known equivalence of $(1+1)$-TQFTs and Frobenius algebras, it also addresses a dual: what conditions on a $(1+1)$-TQFT allow for it to be extended to a $(2+1)$-TQFT? 
		\end{Rem}

		\begin{Rem}\label{rem:signs}
			The proof of our main theorem relies on the functoriality of Bar-Natan's construction.  This holds only over $\bz/2\bz$ -- over $\bz$, the map assigned to a cobordism is only well-defined up to sign.  Indeed, Jacobsson showed that this problem is built in to the definition of Khovanov homology.  So our restriction to $\bz/2\bz$ is out of necessity rather than just convenience.  It is possible that one of the various methods of fixing the sign issues in Khovanov homology will fit naturally into the Khovanov-Floer picture.
		\end{Rem}
\section{Bar-Natan's formulation of link homology}\label{sec:cones}

	This section is a review of the first few sections of \cite{BarNatan2005}.

\begin{figure}
		\begin{centering}
			\includegraphics[width=.5\linewidth]{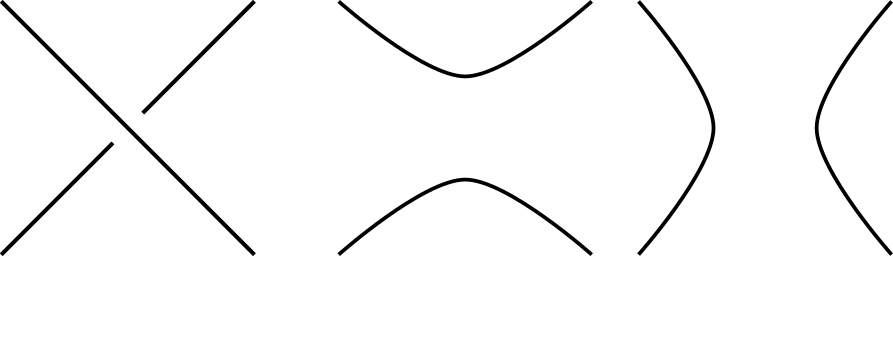}
			\caption{A crossing, its zero-resolution, and its one-resolution.}\label{fig:resolutions}
		\end{centering}
	\end{figure}

	Let $\diagram$ be a link diagram with $c$ ordered crossings.  An element $I$ of $\{0,1\}^c$ is called a \emph{resolution} of $\diagram$.  From a resolution one can produce a new diagram, $\diagram(I)$, also called a resolution, by changing each crossing according to the rule in Figure \ref{fig:resolutions}.  The result is always a diagram without crossings.

	Suppose that $I$ and $I'$ are both resolutions of $\diagram$ which differ only at one entry at which $I$ has $0$ and $I'$ has $1$.  Then there is a cobordism $\Sigma_{I,I'} \co \diagram(I) \to \diagram(I')$ which looks like a saddle near the changed piece and like the identity cobordism everywhere else.  We may describe this saddle (up to isotopy) by drawing an arc (unique up to planar isotopy) on the $0$-resolution of diagram connecting the two sides of the once-crossing.  For a crossing $c$ this arc is called $\ga_c$.  In terms of diagrammatic handle attachments, $\Sigma_{I,I'}$ is a one-handle attachment along $\ga_c$.

	Topologists often describe Khovanov homology as follows: form the \emph{cube of resolutions} of $\diagram$ by drawing the cube $\{0,1\}^c$ and placing the diagram $\diagram(I)$ at the vertex $I$.  Label the edge from $I$ to $I'$ with the cobordism $\Sigma_{I,I'}$.  Apply some recipe to turn each resolution into a module and each cobordism into a linear map.  Bar-Natan formalized the first two steps, and this formalization will form the foundation for our investigation of conic strong Khovanov-Floer theories.

	There is a category $\tildecob$, the \emph{additive closure of $\cob^3$}, whose objects are formal direct sums of objects in $\cob^3$ and whose morphisms are matrices whose entries are formal linear combinations of morphisms in $\cob^3$.   Composition of morphisms is given by matrix multiplication.  Note that $\tildecob$ includes a zero object and zero maps even though $\cob^3$ does not.

	The cube of resolutions of a link diagram $\diagram$ may be understood as an object
	\[
		\BN{\diagram} = \bigoplus_{I \in \{0,1\}^c} \diagram(I)
	\]
	of $\tildecob$ equipped with an endomorphism
	\[
		d = \bigoplus_{I,I'} \Sigma_{I,I'}.
	\]
	It is not hard to show that $d \circ d = 0$ over$\bz/2\bz$.

	\begin{Def}
		Let $\mathcal{C}$ be an additive category.  A \emph{chain complex in $\mathcal{C}$} is an object $C$ along with an endomorphism $d$ of $C$ so that $d \circ d = 0$.  Let $(C,d)$ and $(C',d')$ be complexes.  A morphism $f \in \Hom_{\mathcal{C}}(C,C')$ is called a \emph{chain map} if $f \circ d = d' \circ f$.  Two chain maps $f, g \in \Hom_{\mathcal{C}}(C,C')$ are \emph{chain homotopic} if there is a morphism $H \in \Hom_{\mathcal{C}}(C,C')$ so that 
		\[
			f + g = d' \circ H + H \circ d.
		\]
		$H$ is called a \emph{chain homotopy}.  $f$ is a \emph{chain homotopy equivalence} if $f$ is homotopic to the identity morphism.

		The \emph{homotopy category of complexes over $\mathcal{C}$}, denoted $K(\mathcal{C})$, is the category whose objects are chain complexes over $\mathcal{C}$ and whose morphisms are homotopy classes of chain maps.
	\end{Def}

	To get a manageable category, we need to cut down $\tildecob$ with some relations. Every morphism in $\tildecob$ is a sum of morphisms in $\cob^3$, so we can define relations on $\tildecob$ by defining relations on $\cob^3$.

\begin{figure}
		\begin{centering}
			\includegraphics[width=.75\linewidth]{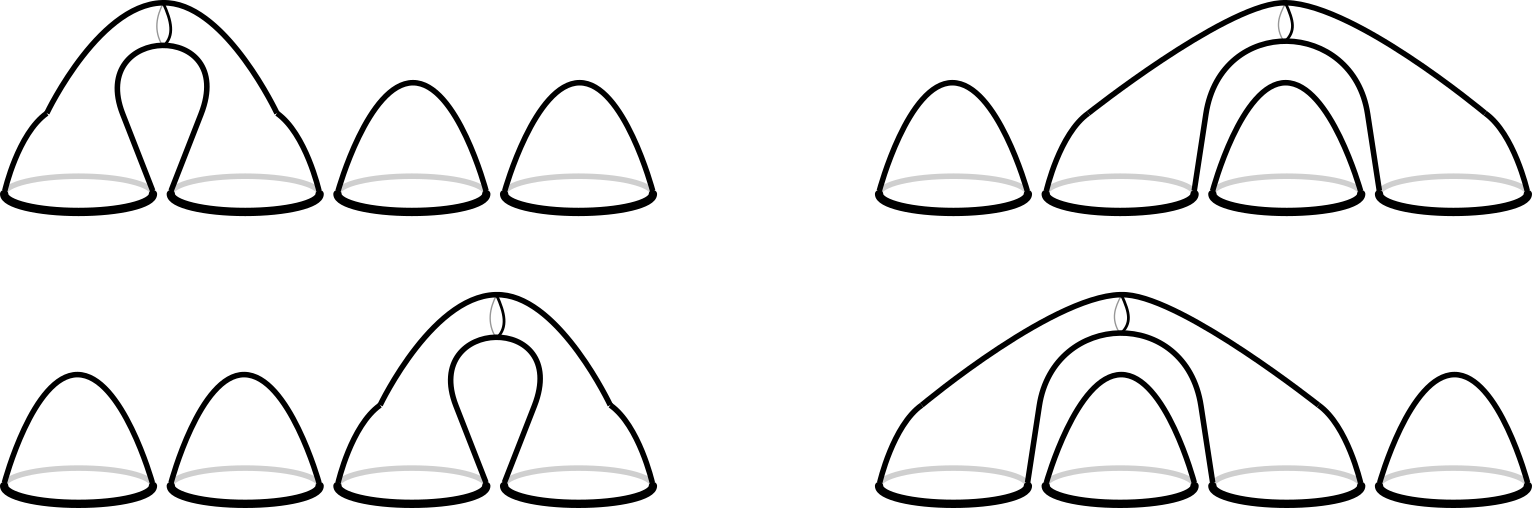}
			\caption{Clockwise from top-left: $\Sigma_{12}$, $\Sigma_{24}$, $\Sigma_{13}$, and $\Sigma_{34}$.  Note that the four disks need not be arranged as in this schematic picture.  All that is required is that they constitute the intersection of some ball with the cobordism $\Sigma$. }
			\label{fig:4tu}
		\end{centering}
	\end{figure}

	\begin{Def}\label{def:relations}
		$\mat(\cob^3)$ is the category whose objects and morphisms are those of $\tildecob$ modulo the following relations on morphisms.
		\begin{description}
			\item[\textbf{The $S$ relation}] any cobordism which includes a closed sphere is equal to zero.
			\item[\textbf{The $T$ relation}] any cobordism which includes a closed torus is equal to zero.
			\item[\textbf{The $4Tu$ relation}] let $\Sigma$ be a cobordism.  Suppose that the intersection of $\Sigma$ with some ball is given by four disks.  Number the disks from $1$ to $4$.  For $i \neq j$ let $\Sigma_{ij}$ be the cobordism given by attaching an unknotted tube between disk $i$ and disk $j$, see \ref{fig:4tu}.  Then $\Sigma_{12} + \Sigma_{34} = \Sigma_{13} + \Sigma_{24}$.
		\end{description}
	\end{Def}

	Every object of $\mat(\cob^3)$ can be written as a sum of objects in $\cob^3$ and every morphism can be written as a sum of cobordisms (up to equivalence) between these objects.  We call these summands \emph{basic}.  $\matcob$ is additive, and the definitions of $\BN{\cdot}$ and $d$ above descend to $\matcob$. 

	\begin{Thm}\label{thm:bn}[Bar-Natan]
		\begin{itemize}
		\item Suppose that $\diagram$ and $\diagram'$ differ by a Reidemeister move.  Then $\BN{\diagram}$ and $\BN{\diagram'}$ are homotopy equivalent by a map $\rho$ which depends on the Reidemeister move.  The isomorphism class of $\BN{\diagram}$ as an object of $K(\matcob)$ is a link invariant.  
		\item 		Let $\Sigma \co \diagram \to \diagram'$ be a diagrammatic cobordism.  It induces a map $F_{\Sigma} \co \BN{\diagram} \to \BN{\diagram'}$.  If $\Sigma'$ is an equivalent diagrammatic cobordism, then $F_{\Sigma} \simeq F_{\Sigma'}$.  In other words, $\BN{\cdot}$ is a functor from $\diag$ to $K(\matcob)$.
		\end{itemize}
	\end{Thm}

	We will describe all these maps in Section \ref{subsec:reidemiester}.

\section{Cobordisms, mapping cones, and functoriality}\label{sec:cobordismsI}

	\subsection{Handle attachment maps and mapping cones}\label{subsec:BNhandleMaps}

		Let $h \in \Hom_{\diag}(\diagram, \diagram')$ be a diagrammatic handle attachment.  Observe that $\diagram$ and $\diagram'$ have the same number of crossings -- in fact, since the handle attachment is supported in the complement of a neighborhood of each crossing, the crossings and resolutions of $\diagram$ and $\diagram'$ are naturally identified.

		Let $\diagram$ be a link diagram with $k$ crossings and $k > 0$.  For each resolution $I$ of $\diagram$, $h$ induces a map $\handle(I) \in \Hom_{\Planar}(\diagram(I), \diagram'(I))$ or equivalently $\handle(I) \in \Hom_{\cob^3}(\diagram(I),\diagram'(I))$.  Define
		\[
			\handle = \bigoplus_{I \in \{0,1\}^k} \handle(I) \in \Hom_{\matcob}(\diagram, \diagram').
		\]
		These are exactly the handle attachment maps defined in \cite{BarNatan2005}.

		The rest of this section is devoted to reframing all of this in terms of mapping cones.  Let $(C,d)$ and $(C',d')$ be chain complexes in $\matcob$.  Let $f \in \Hom(C,C')$ be a chain map.  The \emph{mapping cone of $f$}, denoted $\cone(f)$, is the complex 
		\[
			\left(C \oplus C', \begin{pmatrix} d & f \\ 0 & d'\end{pmatrix}\right)
		\]
		If $\diagram$ has no crossings, then set
		\[
			\BN{\diagram}' = (\diagram, 0).
		\] 
		If $\diagram$ has $k$ crossings, label (and thus order) them by $\{c_1, \ldots, c_k\}$.  Let $\diagram_{k-1,0}$ and $\diagram_{k-1,1}$ be the diagrams given by resolving the crossing $c_k$.  Let $\ga_{c_k}$ be the arc assigned to the crossing $c_k$ as in the previous section. Write $\handle[k]$ for $\handle[\ga_{c_k}]$.  Define $\BN{\diagram}' = \cone(\handle[k])$.  We call this the \emph{iterated mapping cone} description of $\BN{\diagram}$.  The following proposition is well-known.

		\begin{Prop}\label{prop:coneToCube1}
			$\BN{\diagram} = \BN{\diagram}'$.  In other words, the cube of resolutions and iterated mapping cone descriptions agree.
		\end{Prop}

		\begin{proof}
			For diagrams with no crossings, the descriptions are identical.  Suppose that the proposition holds for diagrams up to $k-1$ crossings.  Let $\diagram$ be a $k$-crossing diagram with crossings $\{c_1, \ldots, c_k\}$.  Then $\BN{\diagram}' = \cone(\handle[k] \co \BN{\diagram_{k-1,0}}' \to \BN{\diagram_{k-1,1}}')$.  By hypothesis, $\BN{\diagram_{k-1,0}}'$ and $\BN{\diagram_{k-1,1}}'$ are both equivalent to cubes of resolutions:
			\begin{align*}
				\BN{\diagram_{k-1,0}}' &= \left(\bigoplus_{I \in \{0,1\}^{k-1}} \diagram_{k-1,0}(I), \partial_{0,k-1}\right) \\
				\BN{\diagram_{k-1,1}}' &= \left(\bigoplus_{I \in \{0,1\}^{k-1}} \diagram_{k-1,1}(I), \partial_{1,k-1}\right)
			\end{align*}
			So 
			\[
				\BN{\diagram}' = \bigoplus_{I \in \{0,1\}^k} \diagram(I)
			\]
			with differential given by $d' = d'_{0,k-1} \oplus d'_{1,k-1} \oplus \handle[k]$.  Now it is clear that $d'$ agrees with $d$.
		\end{proof}

		Observe that the mapping cone description does not use on the fact that $\handle[k]$ is a sum of maps between $\diagram_{k-1,0}(I)$ and $\diagram_{k-1,1}(I)$.  We could write $\handle[k](I)$, the restriction of $\handle[k]$ to $\diagram_{k-1,0}(I)$, as a map
		\[
			\handle[k](I) \co \BN{\diagram_{k-1,0}(I)} \to \BN{\diagram_{k-1,1}}.
		\]
		Let $\handle$ be a diagrammatic handle attachment.  Write $\BN{\diagram} = \cone(\handle[k] \co \diagram_0 \to \diagram_1)$.  The arc $\gamma_{c_k}$ is disjoint from the support of $\handle$, so $\handle[k]$ and $\handle$ commute in $\Planar$ and therefore in $\cob^3$.  It follows that
		\[
			\handle \co \cone(\handle[k] \co \BN{\diagram_0} \to \BN{\diagram_1}) = \cone(\handle[k] \co (\BN{\handle\diagram_0}) \to (\BN{\handle\diagram_0})).
		\]
		\begin{figure}
			\begin{centering}
				\includegraphics[width=.75\linewidth]{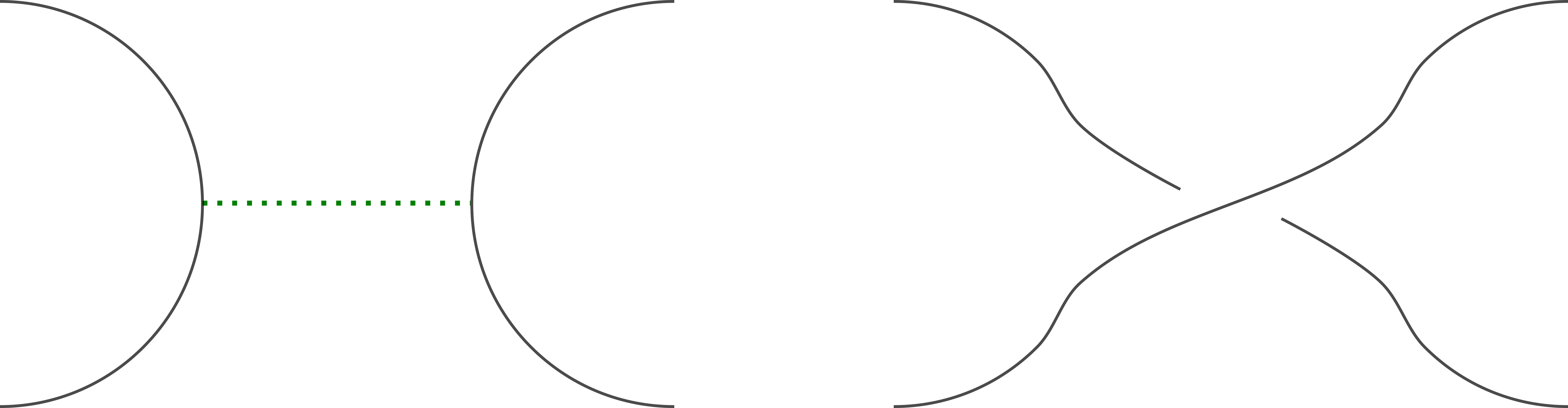}
				\caption{On the left, a link diagram in black and the curve $\gamma$ in green.  On the right, $\diagram_\gamma$.}
				\label{fig:addCrossing}
			\end{centering}
		\end{figure}

	\subsection{Mapping cones and strong Khovanov-Floer theories}\label{subsec:cones}

		Recall that a strong Khovanov-Floer theory is conic if the complex it assigns to a link diagram is the mapping cone of the one-handle attachment map given by any crossing.  In other words, the total complex $\hkft{\diagram}$ is an iterated mapping cone.

		\begin{Prop}\label{prop:coneToCube2}
			   Let $\hkft$ be a conic strong Khovanov-Floer theory.  Then $\hkft$ admits a cube of resolution descriptions.  That is, for any link diagram $\diagram$, the group $\hkft(\diagram)$ is
			   \[
			   	\bigoplus_{I} \hkft(\diagram(I))
			   \]
			   and the differential on $\hkft(\diagram)$ is a sum of maps $\partial_{I,J}$ with $I \leq J$.
		\end{Prop}

		\begin{proof}
			Let $\diagram$ be a link diagram with crossings $\{c_1, \ldots, c_k\}$.  Write $f^{\hkft}_{c_k} = \hkft(\handle[c_k])$.

			This follows from the proof of Proposition \ref{prop:coneToCube1}with two alterations.  The first is that, for any resolution $I$, the complex $\hkft(\diagram(I))$ may have a differential.  The second difference is that it may not be possible to write $f_{c_k}$ as a sum of maps $f^\hkft_{c_k}(I) \co \hkft(\diagram_0(I)) \to \hkft(\diagram_1(I))$.  Rather, one has
			\begin{align*}
				f^\hkft_{c_k}(I, J) &\co \hkft(\diagram_0(I)) \to \hkft(\diagram_1(J)) \\
				f^\hkft_{c_k}(I) &= \sum_{I \leq J} f^\hkft_{c_k}(I,J).
			\end{align*}
			In other words, there may be `edges' in the cube of length zero and length greater than one.  Nevertheless, the observation following the proof of Proposition \ref{prop:coneToCube1} makes clear that the proof carries through without alteration.
		\end{proof}

		So if $\hkft$ is conic, then $\hkft(\diagram)$ has a homological filtration: if $x \in \hkft(\diagram(I))$, then the homological grading of $x$ is the sum of the entries in $I$.  (This grading may need to be shifted to be a link invariant -- in any case, there is a well-defined relative grading.)

		Now we can connect strong Khovanov-Floer theories to Khovanov-Floer theories.  The essential point is that a cube of resolutions description allows us to easily recognize associated graded objects.

		\begin{Prop}\label{prop:homotopyToRegular}
			Let $\hkft$ be a conic strong Khovanov-Floer theory with underlying Frobenius algebra $R$. It defines a Khovanov-Floer theory: to a link diagram $\diagram$ it assigns the spectral sequence given by the homological filtration on $\hkft(\diagram)$.  To a one-handle attachment it assigns the map of spectral sequences induced by $\hkft$.
		\end{Prop}
		Note that the Lee-Bar-Natan spectral sequence is not induced by the homological filtration.
		\begin{proof}
			Write $\partial$ for the differential on $\hkft(\diagram)$.  Write $\partial_i$ for the component of $d$ with homological degree $i$.  Consider
			\[
				H(\hkft(\diagram),\partial_0) = \bigoplus_{I \in \{0,1\}^k} H(\diagram(I),\partial_0(I))
			\]
			We work by induction on $k$ to show that $(\partial_1)_*$ is equal to the differential on $R$ link homology.  Condition \ref{kft:frobenius} in the definition of strong Khovanov-Floer theory implies the claim for diagrams with one crossing.  For a diagram $\diagram$ with $k$ crossings, $k>1$, write $\hkft(\diagram)$ as the cone of $f_{c_k}$.  Then $H(\hkft(\diagram),d_0)$ is a cone on $(f_{c_k})_*$, the map induced by $f_{c_k}$.  It suffices to show that the length one part of $(f_{c_k})_*$ agrees with the map given by $R$.

			Let $I$ and $I'$ be two resolutions which differ only at the $k$th crossing.  By hypothesis, $\diagram$ has at least two crossings, so there is another crossing (say, $c_1$) at which $I$ and $I'$ agree.  Write $H(\hkft(\diagram),d_0)$ as the cone on $(f_{c_1})_*$.  The component of $(f_{c_k})_*$ from $H(\hkft(\diagram(I)),d_0)$ to $H(\hkft(\diagram(I')),d_0)$ must agree with the $R$-link homology map by the induction hypothesis. 
		\end{proof}
	\subsection{Strong Khovanov-Floer theories and planar diagrams}\label{subsec:strongKFandPlanar}

		It's time to connect strong Khovanov-Floer theories to Bar-Natan's category.

		\begin{Prop}\label{prop:skftMeetsBN}
			Every strong Khovanov-Floer theory satisfies the $S$, $T$, and $4Tu$ relations.  That is, if $\hkft$ is a strong Khovanov-Floer theory and $\Sigma$ is a diagrammatic cobordism containing a sphere or torus, then $\hkft(\Sigma) = 0$.  Using the notation of Definition \ref{def:relations}, 
			\[
				\hkft(\Sigma_{13}) + \hkft(\Sigma_{24}) = \hkft(\Sigma_{12}) + \hkft(\Sigma_{34}).
			\]
		\end{Prop}

		\begin{proof}
			Suppose that $\Sigma$ contains a sphere.  Arrange that the sphere has the standard handle decomposition with just a zero- and two-handle using handleswaps and movie move 15.  Conditions \ref{kft:crossingless} and \ref{kft:handleKunneth} of Definition \ref{def:gkft} imply that $\hkft(\Sigma) = 0$.  The same argument applies to the $T$ relation.

			For the $4Tu$ relation, see the proof of Proposition 7.2 in \cite{BarNatan2005}.  Again, conditions \ref{kft:crossingless} and \ref{kft:handleKunneth} along with a short calculation verify the relation.
		\end{proof}

		\begin{Rem}\label{rem:thanksGordana}
			It is interesting to consider the four-dimensional geometry of Bar-Natan's relations.  The $S$ relation is familiar: let $W \co Y \to Y'$ be a four-dimensional cobordism so that there is some $h \in H_2(W)$ so that $h$ can be represented by an embedded sphere of self-intersection number zero.  Then the map assigned to this cobordism by any flavor of instanton Floer homology is zero \cite{Kronheimer1995}.  If $\Sigma$ is a link cobordism which contains a sphere, then its branched double cover will contain a sphere of self-intersection zero.  (I am thankful to Gordana Mati{\`c} for pointing this out.)

			The twice-punctured torus appears in Section 8.3 of \cite{Kronheimer2011}.  It is essential in the identification of cobordism maps on unlinks with the maps in Khovanov homology in the following section.  A four-dimensional interpretation of the $4Tu$ relation is less clear to us.
		\end{Rem}

		Here is some notation necessary for our main technical result.  Let $\Kom$ for the homotopy category of complexes with coefficients in $\bF$.  Let $\diag_\circ$ be the subcategory of $\diag$ with the same objects but with morphisms restricted to movies without Reidemeister moves.  Write $T(\circ)$ for the smallest full subcategory of $\matcob$ whose objects include every finite configuration of circles in $\br^2$ and so that if $\handle \in \Hom(O_1, O_2)$ is a one-handle attachment map between two objects of $T(\circ)$, then $\cone(\handle)$ is an object of $T(\circ)$.

		\begin{Prop}\label{prop:mainTechnical}
			\begin{enumerate}	
				\item\label{i:tech1} A conic strong Khovanov-Floer theory $\hkft$ defines a functor from $\diag_\circ$ to the category of filtered complexes with coefficients over $\bF$.  The functor sends an object $\diagram$ to $\hkft(\diagram)$.  Let $\Sigma \co \Hom_{\diag_\circ}(\diagram, \diagram')$.  Write it as a composition of isotopies and diagrammatic handle attachments
				\[
					\Sigma = \Sigma_n \circ \cdots \circ \Sigma_1.
				\]
				Then
				\[
					\hkft(\Sigma) = \hkft(\Sigma_n) \circ \cdots \circ \hkft(\Sigma_1).
				\]
				\item\label{i:tech2} The functor extends linearly to a functor $\mathcal{F}^{\circ}_{\hkft}$ from $\matcob$ to $\Kom$.
				\item\label{i:tech3} Let $\partial$ be the differential on $\hkft(\diagram)$ and let $d$ be the differential on $\BN{\diagram}$.  Let $\partial_0$ be the part of $\partial$ with homological degree $0$.  Then $\mathcal{F}^{\circ}_{\hkft}(d) = \partial^0 + \partial$.
				\item\label{i:tech4} The functor $\mathcal{F}^{\circ}_{\hkft}$ induces a functor $\hkft$ from $T(\circ)$ to $\Kom$ so that $\hkft = \mathcal{F}^{\circ}_{\hkft} \circ \BN{-}$.  This notation makes sense: the object assigned to $\diagram$ by the functor $\hkft$ is the complex $\hkft(\diagram)$, and similarly for morphisms.
			\end{enumerate}
		\end{Prop}

		\begin{proof}
			\begin{enumerate}
				\item Maps for handle attachments and planar isotopies are part of the data of a strong Khovanov-Floer theory.  Suppose that 
				\[
					\Sigma = \Sigma'_m \circ \cdots \circ \Sigma'_1
				\]
				for some other elementary cobordisms $\Sigma'_i$.  Then by Lemma \ref{lem:movies2}, the movies given by these decompositions differ by some sequence of distant handleswaps and applications of movie move 15.  By Condition \ref{def:relations} of Definition \ref{def:gkft}, any application of these moves preserves the chain homotopy class of the induced map.  So $\hkft(\Sigma)$ is well-defined.  

				Compatability with composition and identity maps is clear.
				
				\item Linearly extend the functor above to $\tildecob$.  Proposition \ref{prop:skftMeetsBN} shows that it descends to $\matcob$.  
				
				\item For diagrams without crossings $\mathcal{F}^{\circ}_{\hkft}(d) = 0$ and $\partial^0 = \partial$.  For a diagram with $k > 0$ crossings, $\BN{\diagram}$ is the cone of the map 
				\[
					\handle[k] \co \BN{\diagram_0} \to \BN{\diagram_1}.
				\]  
				Write $d_i$ for the differential on $\BN{\diagram_i}$ and $\partial_i$ for the differential on $\hkft(\diagram_i)$, $i = 0, 1$.  Therefore $d = d_0 \oplus d_1 \oplus \handle[k]$.  By hypothesis, $\hkft(d_0)$ and $\hkft(d_1)$ agree with $\partial_0 + \partial_0^0$ and $\partial_1 + \partial_1^0$, and by conicity $\partial = \partial_0 + \partial_1 + \hkft(\handle[k])$.  Observe that $\partial^0 = \partial_1^0 + \partial_0^1$.  It follows that $\hkft(d) = \partial + \partial^0$.

				\item The definition on objects is clear.  The real content of this statement is that if $\Sigma_1, \Sigma_2 \in \Hom_{T(\circ)}(O,O')$ and $\BN{\Sigma_1} \simeq \BN{\Sigma}$, then $\hkft(\Sigma_1) = \hkft(\Sigma_2)$.  It suffices to prove this in the case that $O$ and $O'$ are finite collections of circles -- the result for mapping cones of one-handle attachments follows as in the last item.

				$\Sigma_1$ and $\Sigma_2$ are equal in the quotient of $\cob^3$ by the $S$, $T$, and $4Tu$ relations.  This implies that, after adding some spheres and tori to $\Sigma_1$ and $\Sigma_2$, they are isotopic in $\br^2 \times I$ after some applications of the $4Tu$ relation.  Call these isotopic cobordisms $\Sigma'_1$ and $\Sigma'_2$. Lemma \ref{lem:movies2} implies that any two movies $M(\Sigma'_1)$ and $M(\Sigma'_2)$ are related by distant handleswaps and movie move 15.  So Proposition \ref{prop:skftMeetsBN} and Condition \ref{def:relations} imply that $\hkft(\Sigma_1) \simeq \hkft(\Sigma_2)$.
			\end{enumerate}	
		\end{proof}

	\subsection{Reidemeister maps and invariance}\label{subsec:reidemiester}

		\begin{figure}
			\centering
			\huge
			\def\svgwidth{.75\linewidth}{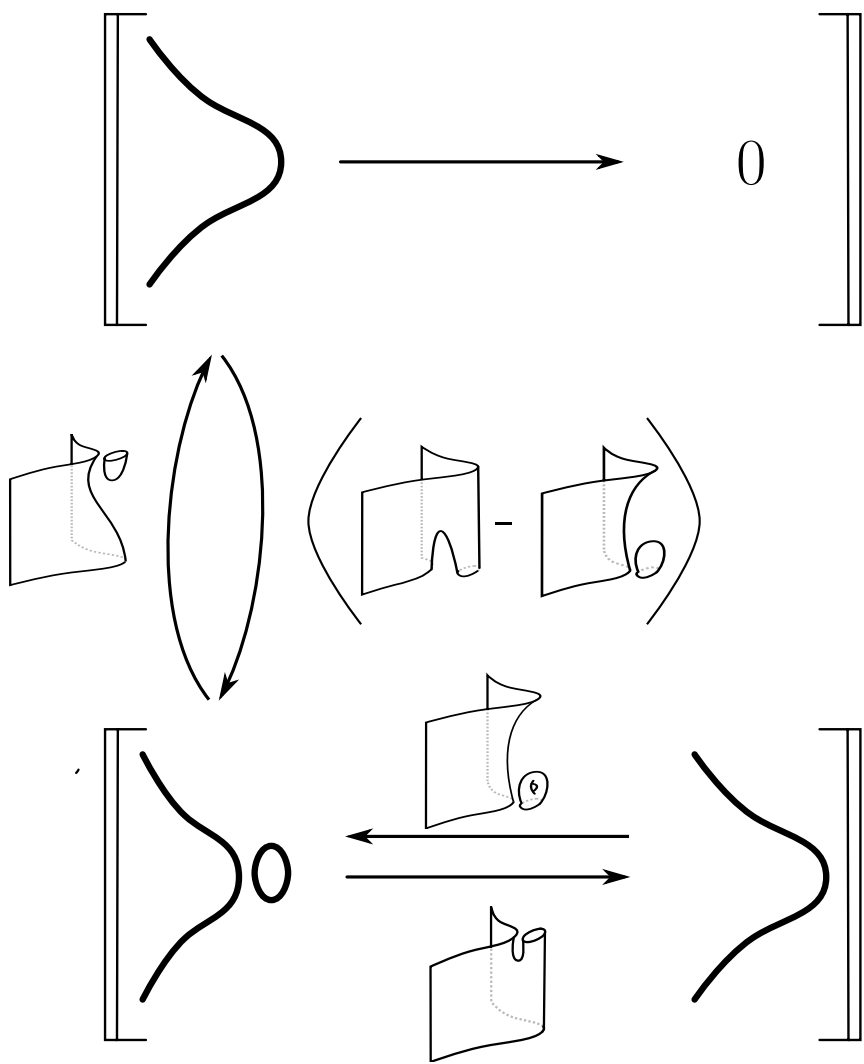}
			\caption{Reidemeister 1 invariance for $\BN{\cdot}$.  The vertical maps are chain maps.  The right-to-left horizontal map on the bottom is a chain homotopy.}
			\label{fig:BNR1}
		\end{figure}

		Let $\diagram$ be a link diagram with $k$ crossings and let $\diagram'$ be the result of applying a left-handed Reidemeister 1 move to $\diagram$.  Figure \ref{fig:BNR1} shows a map $\BN{\diagram} \to \BN{\diagram'}$.  There is an obvious identification of $k$ of the crossings of $\diagram'$ with those of $\diagram$.  So 
		\[
			\BN{\diagram'} = \cone(\handle[k+1] \co \diagram'_0 \to \diagram'_1)
		\]
		where $\diagram'_0$ and $\diagram'_1$ are the diagrams obtained by taking the $0$- and $1$-resolution of $c_{k+1}$, respectively.  This description is local in the sense that all the diagrams are identical outside of the support of the Reidemeister move.  The cobordisms in the figure are understood to be extended to the rest of the diagram by the identity.  Bar-Natan shows that the vertical maps are chain maps and that they form a chain homotopy equivalence.  (The right-to-left horizontal map is the homotopy.)  He produces similar maps for Reidemeister 2 and 3 moves.\footnote{One can also prove Reidemeister 3 invariance from a comparison lemma for mapping cones (Lemma 4.5 of \cite{BarNatan2005}). This lemma holds in any additive category and therefore in $\Kom$.}

		Proposition \ref{prop:mainTechnical} implies that each individual cobordism defines a chain map in a strong Khovanov-Floer theory.  Bar-Natan's proof that the ensemble defines a chain map $\BN{\diagram} \to \BN{\diagram'}$ carries through without alteration.  Call the resulting map $\hkft(\rho)$ and call the reverse map $\hkft(\rho')$.  Then
		\[
			\rho' \circ \rho = \Id + dH + Hd
		\]
		and therefore
		\begin{align*}
			\hkft(\rho') \circ \hkft(\rho) &= \hkft(\Id) + \hkft(dH) + \hkft(Hd)
			\\
			&= \Id + \partial\hkft(H) + \hkft(H)\partial + \partial^0\hkft(H) + \hkft(H)\partial^0
		\end{align*}
		The critical observation is that \emph{the homotopies $H$ consist only of zero- and two-handle attachments}.  These must commute with $\partial^0$ because they define chain maps with homological degree zero.

		\begin{Prop}\label{prop:reidemeister}
			Let $\hkft$ be a conic strong Khovanov-Floer theory.  Let $\diagram$ and $\diagram'$ be link diagrams which differ by some sequence of Reidemeister moves.  Then there is a filtered homotopy equivalence $\rho \co \hkft(\diagram) \to \hkft(\diagram')$.  This map induces a map on the associated graded which agrees with the usual map on the Khovanov chain complex.
		\end{Prop}

		\begin{proof}
			We have already proved everything except the last sentence, which follows by an argument similar to that of \ref{prop:homotopyToRegular}.
		\end{proof}

		\begin{Rem}
			Bar-Natan develops a planar algebra theory under which these pictures can be taken literally rather than locally; they are objects in (the additive closure of) a category of tangles.  This puts the structure of a planar algebra on Khovanov homology and its variants.  We are not aware of planar algebra structures in any other Khovanov-Floer theory.  As a middle ground, there are versions of many Khovanov-Floer theories for two-sided tangles (e.g. braids).  The author and John Baldwin have developed an invariant of two-sided tangles (e.g. braids) for \Szabo's homology theory, \cite{Baldwin2018}.  This suggests that strong Khovanov-Floer theories should all admit some planar algebra structure.  We will not pursue this angle further here.
		\end{Rem}

	\subsection{Functoriality}

		Let $M \in \Hom_{\diag}(\diagram, \diagram')$ be a movie.  Write $M$ as a composition of elementary diagrammatic cobordisms, i.e. Reidemeister moves, handle attachments, and isotopies:
		\[
			M = M_n \circ \cdots \circ M_1.
		\]
		Define
		\[
			\BN{M} = \BN{M_n} \circ \cdots \circ \BN{M_1}.
		\]
		Let $I$ be a resolution of $\diagram$.  Let $\BN{M}(I)$ be the restriction of $\BN{M}$ to $\BN{\diagram(I)}$.  Observe that there is a collection of movies $\{M^1(I), \ldots, M^m(I)\}$ \emph{without Reidemeister moves} so that 
			\[
				\sum_{j=1}^m \BN{M^j(I)} = \BN{M}(I).
			\]
		If $M$ has no Reidemeister moves then the collection is just $\{M\}$ and there is nothing to prove.  If $M$ has Reidemeister moves, then replace the Reidemeister maps from the previous section.  These maps are each sums of restrictions of diagrammatic handle attachments.  Let $\mathcal{K}$ be a conic, strong Khovanov-Floer theory.  Define
		\[
			\hkft(M_i) = \mathcal{F}^\circ_{\hkft}(\BN{M_i})
		\]
		and
		\[
			\hkft(M) = \hkft(M_n) \circ \cdots \circ \hkft(M_1).
		\]
		To be concrete: if $M_i$ is a diagrammatic handle attachment, then $\hkft(M_i)$ is determined by the definition of $\hkft$.  If $M_i$ is a Reidemeister move, then $\hkft(M_i)$ is determined by Proposition \ref{subsec:reidemiester}. 
		\begin{Thm}\label{thm:main}
			Every conic strong Khovanov-Floer theory is functorial: if $M$ and $M'$ are movies representing isotopic link cobordisms, then $\hkft(M)$ is filtered chain homotopic to $\hkft(M')$.
		\end{Thm}
		\begin{proof}
			Suppose that $M$ and $M'$ are movies representing isotopic cobordisms from $\diagram$ to $\diagram$.  Theorem \ref{thm:bn} shows that 
			\begin{equation}\label{eqn:1}
				\BN{M} + \BN{M'} = d' \circ H + H \circ d
			\end{equation}
			where $d'$ is the differential on $\BN{\diagram'}$ and $H \in \Hom_{\diag^\circ}(\diagram,\diagram')$.  Write $H = \sum_{i=1}^n \BN{h_i}$ where the $h_i$ are restrictions of elementary diagrammatic cobordisms.  Let $\partial$ and $\partial'$ be the differentials on $\hkft(\diagram)$ and $\hkft(\diagram')$.  Proposition \ref{prop:mainTechnical} shows that $\hkft(d) = \partial$ and $\hkft(d') = \partial'$.  We have defined $\hkft(M) = \mathcal{F}^\circ_{\hkft}(\BN{M})$.  So applying $\mathcal{F}^\circ_{\hkft}$ to both sides of the \ref{eqn:1} yields
			\[
				\hkft(M) + \hkft(M') = \partial'\hkft(H) +  \hkft(H)\partial + (\partial'^0\hkft(H) + \hkft(H)\partial'^0)
			\]
			So $\hkft(M)$ and $\hkft(M')$ are chain homotopic via $H = \hkft(d'),$ a filtered chain homotopy, as long as $\hkft(H)$ intertwines $\partial'^0$ and $\partial^0$.

			The proof of Theorem \ref{thm:bn} sorts the movie moves into three types.  The Type I moves say that a cobordism consisting of a Reidemeister move and its reverse is isotopic to a cylinder.  The Type II moves have on one side a complicated sequence of Reidemeister moves and on the other a cylinder.  The Type III moves involve handle attachments.  If $M$ and $M'$ are related by a move of Type I, then we have already shown that $\hkft(H)$ intertwines $\partial'^0$ and $\partial^0$ in Proposition \ref{prop:reidemeister}.  With the exception of movie move 15, the cobordisms on each side of a Type III move are actually equal modulo the $S$, $T$, and $4Tu$ relations -- the homotopy $H$ is zero.  If $M$ and $M'$ are related by movie move 15, then $H$ is essentially the Reidemiester 3 homotopy.

			That leaves the movie moves of Type II.  The key lemma is this: let $t \subset \diagram$ be a tangle without crossings or closed components in a link diagram.  Let $\Sigma \co \diagram \to \diagram$ be a diagrammatic cobordism without handle attachments whose support is in $t$.  Then $\BN{\Sigma} = \Id$ -- the proof uses only the $S$ and $T$ relations, and it therefore holds for $\hkft(\Sigma)$.  

			Fix some Type II move and call its complicated side $M_L$.  Its initial frame may be viewed as a braid $\beta$ in some link diagram $\diagram$.  There is an isotopy consisting only of Reidemeister 2 moves from a tangle $t$ without crossings or closed components to $\beta\beta^{-1}$.  Now perform the sequence of diagrammatic cobordisms prescribed by $M_L$, then undo $\beta\beta^{-1}$ back to $t$.  This whole sequence must induce the identity map on $\hkft(\diagram)$, i.e.
			\[
				\rho' \circ \BN{M_L} \circ \rho = \Id.
			\]
			It follows that
			\[
				\BN{M_L} \simeq \rho' \circ \rho \simeq \Id
			\]
			where each homotopy is given by Reidemeister 2 homotopies.  Therefore $\hkft(H)$ intertwines $\partial^0$ and $\partial'^0$ for Type II moves.
		\end{proof}

\section{Link homologies as strong Khovanov-Floer theories}\label{sec:examples}
	
		In this section we show that $\Bn$ and $\CSz$ are conic, strong Khovanov-Floer theories.  We will write $\hkft$ and $\handle[\gamma]$ without specifying the particular theory when it is understood.   

		Every strong Khovanov-Floer theory we consider will be conic.  So the recipe for constructing one-handle attachment maps will be common among all of them: let $\diagram$ be a link diagram and $\gamma$ a planar arc in $\diagram$.  Let $\diagram'$ be the diagram given by adding a crossing along $\gamma$ as in Figure \ref{fig:addCrossing}.  Then $\hkft(\diagram')$ is the cone of some map.  This is the one-handle attachment map $\handle[\gamma]$.  \emph{A prioi} this may still depend on some choices of auxiliary data.

		For \emph{conic} strong Khovanov-Floer theories there is already some redundancy within the definition.  For this section, let a \emph{pretheory} be a rule which satisfies some of the conditions in Definition \ref{def:gkft}.

		\begin{Prop}\label{prop:universalDistanceRule}
			Let $\hkft$ be a conic pretheory which satisfies conditions \ref{kft:coherence} and \ref{kft:handles} of Definition \ref{def:gkft}.  Then it satisfies handleswap invariance.
		\end{Prop}
		\begin{proof}
			Let $\gamma$ and $\gamma'$ be two arcs for one-handle attachments on the link diagram $\diagram$.  Let $\diagram''$ be the link diagram given by adding crossings at $\gamma$ and $\gamma'$ as in \ref{fig:addCrossing}.  By hypothesis, $\hkft(\diagram'')$ is a chain complex.  Write $\diagram''_{ij}$ with $i,j \in \{0,1\}$ for the diagrams given by resolving these two crossings according to $i$ and $j$.  Now consider the component of $\partial^2$ extending from $\hkft(\diagram_{00})$ to $\hkft(\diagram_{11})$.  This can be written as a sum
			\[
				\handle[\gamma'] \circ \handle[\gamma] + \handle[\gamma] \circ \handle[\gamma'] + \partial H + H \partial
			\]
			where $H$ is the component of $\partial$ from $\hkft(\diagram_{00})$ to $\hkft(\diagram_{11})$.  Therefore, absuing notation a bit, $\handle[\gamma] \circ \handle[\gamma'] \simeq \handle[\gamma] \circ \handle[\gamma']$.  
		\end{proof}

\begin{figure}
			\centering
			\def\svgwidth{.75\linewidth}{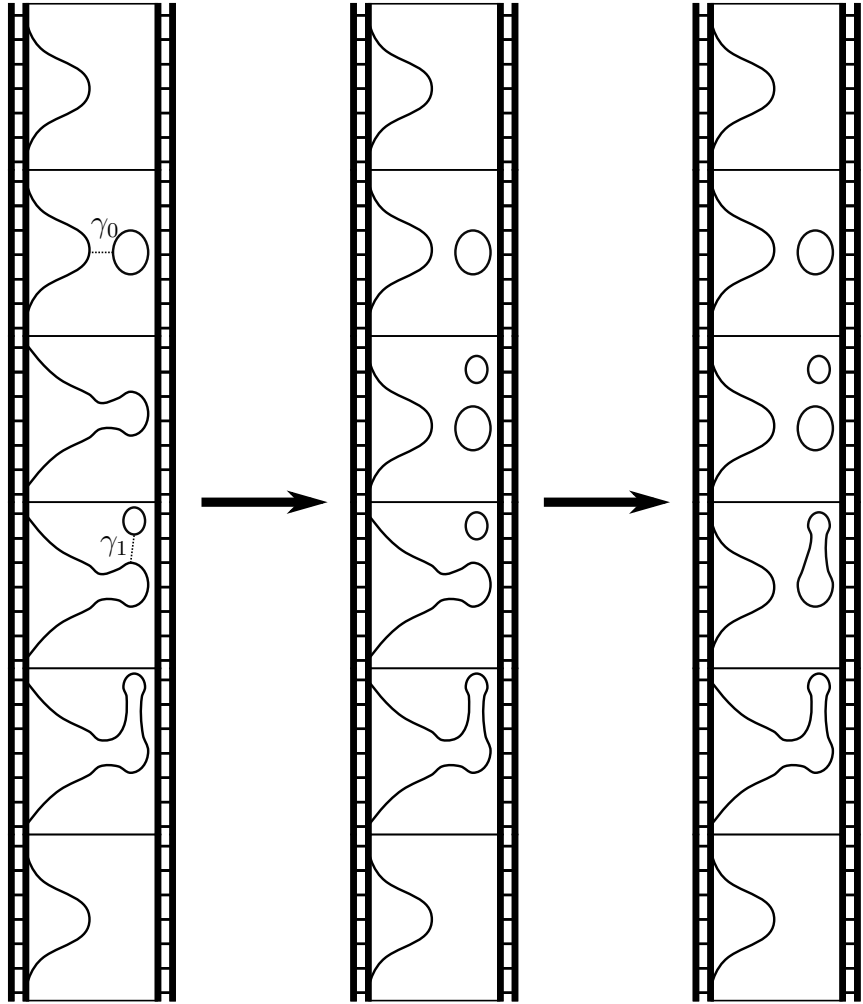}
			\caption{On the left, a movie representing $\handle[\gamma_1]Z_{1}\handle[\gamma_0]  Z_0$.  The middle represents $\handle[\gamma_1]\handle[\gamma_0] Z_1 Z_0$.  The right represents $\handle[\gamma_0]\handle[\gamma_1] Z_{1}Z_0$.}
			\label{fig:double16}
		\end{figure}

		\begin{Prop}\label{prop:conicImplies16}
			Let $\hkft$ be a conic pretheory which satisfies every condition except movie move 15 invariance.  Then $\hkft$ satisfies movie move 15 invariance.
		\end{Prop}
		\begin{proof}
			By conditions \ref{kft:coherence} and \ref{kft:handles}, we are justified in talking about groups and handle attachment maps.  Note that \ref{kft:kunneth} implies that, if $\diagram$ is the empty diagram, then $\hkft(\diagram) \cong \bF$.

			Let $\diagram$ be a link diagram.  Let $\diagram'$ be the result of a zero- or two-handle attachment on $\diagram$.  There is an obvious bijection between the resolutions of $\diagram$ and those of $\diagram'$.  We first show that the handle attachment map $f \co \hkft(\diagram) \to \hkft(\diagram')$ decomposes (up to chain homotopy) as a sum
			\begin{align*}
				f &= \sum_{I} f_I \\
				f_I &\co \hkft(\diagram(I)) \to \hkft(\diagram'(I)).
			\end{align*}
			Say that such a map \emph{decomposes resolution by resolution}.  This follows from Condition \ref{kft:handleKunneth}. Suppose that $f$ is a zero-handle attachment map.  Then $f$ is chain homotopic to the composition
			\[
				\hkft(\diagram) \cong \hkft(\diagram) \map{\Id \otimes f_\circ} \bF \to \hkft(\diagram) \otimes \hkft(\circ) \cong \hkft(\diagram \cup \circ) \cong \hkft(\diagram').
			\]
			If $f$ is a two-handle attachment map, then $f$ is chain homotopic to the reverse composition.  Every map in that composition decomposes resolution by resolution.

			Let $\Sigma$ be a diagrammatic cobordism which is equivalent to the identity by a single application of movie move 15.  The homologically grading-preserving part of this map is the corresponding map in some Frobenius algebra by Condition \ref{kft:frobenius}.  Therefore this map is the identity.  Any filtered map on a bounded filtered complex whose graded part is the identity is a quasi-isomorphism.  (This is a well-known exercise; see Lemma 24 of \cite{Saltz2017} for a solution.)  So $\hkft(\Sigma)$ is a quasi-isomorphism.

			Now we show that $\hkft(\Sigma) \circ \hkft(\Sigma) \simeq \hkft(\Sigma)$.  Then Lemma \ref{lem:lessDumbAlgebra} finishes the proof.  Suppose for now that $\Sigma$ consists of a canceling zero- and one-handle.  Then $\hkft(\Sigma) = \handle[\gamma] \circ Z$ where $Z$ is the zero-handle attachment map.  The key observation is that
			\[
				\hkft(\Sigma)\hkft(\Sigma) =  \handle[\gamma_1]Z_{1}\handle[\gamma_0]Z_0
			\]
			where $\gamma_0$ merges $\circ$ into $\diagram$ and $\gamma_1$ merges $\circ'$ to $\circ$.  See Figure \ref{fig:double16} for a film study.  $Z_1$ commutes with $\handle[\gamma_0]$ because their supports are disjoint, so
			\[
				\hkft(\Sigma)\hkft(\Sigma) \simeq \handle[\gamma_0]\handle[\gamma_1]Z_{1}Z_{0}.
			\]
			Because $\hkft$ satisfies Condition \ref{kft:frobenius},
			\[
				\handle[\gamma_1]Z_{1}Z_{0} \simeq Z_{2}
			\]
			where $Z_{2}$ is a zero-handle whose support intersects the support of both $Z_0$ and $Z_1$.  So
			\[
				\hkft(\Sigma) \circ \hkft(\Sigma) \simeq \handle[\gamma_0] \circ Z_{2}
			\]
			where $\gamma_0$ now attaches $\circ_2$ to $\diagram$.  Therefore $\hkft(\Sigma) \circ \hkft(\Sigma) \simeq \hkft(\Sigma)$.

			For the two-handle case, just turn this argument upside-down: show that 
			\[
				T_{\circ_1}\handle[\gamma_1]T_{\circ_0}\handle[\gamma_0] \simeq T_{\circ_1}T_{\circ_0}\handle[\gamma_1]\handle[\gamma_0]
			\]
			where the $T$s are two-handle attachment maps, and $\gamma_0$ and $\gamma_1$ are arcs which split off a little basepoint-anchored circle with $\gamma_1$ nested inside $\gamma_0$.  Now shuffle around the maps and use the Khovanov merge map to show that $\hkft(\Sigma) \circ \hkft(\Sigma) \simeq \hkft(\Sigma)$.
		\end{proof}

		\begin{Lem}\label{lem:lessDumbAlgebra}
			Let $C$ be a chain complex over $\bF$ and let $f \co C \to C$ be a quasi-isomorphism.  Suppose that $f^2 \simeq f$.  Then $f \simeq \Id$.
		\end{Lem}
		\begin{proof}
			The facts that $(f^2)_* = f_*$ and that $f_*$ is an isomorphism imply that $f_* = \Id_*$.

			Now we show that the projection $\pi \co C \to H(C)$ and the inclusion $\iota \circ H(C) \to C$ form a homotopy equivalence. (This seems to be well-known in some communities; we give a proof for completeness.)  Let $d$ be the differential on $C$.  $\pi$ and $\iota$ are chain maps, and $\pi \circ \iota$ is the identity on $H(C)$.  Working over a field, we have
			\[
				C \cong \ker(d) \oplus C/\ker(d) \cong H(C) \oplus \im(d) \oplus C/\ker(d).
			\]
			Let $h$ be the map which vanishes on $H(C)$ and $C/\ker(d)$ and which acts as the isomorphism $\im(d) \cong C/\ker(d)$ on $\im(d)$.  Notice that the isomorphism $C/\ker(d) \cong \im(d)$ is the restriction of $d$.  Then
			\[
				\iota \circ \pi = \Id + dh + hd.
			\]
			Certainly this holds on $H(C)$.  On $\im(d)$, the left side vanishes, $hd = 0$, and $\Id = dh$.  On $C/\ker(d)$, the left side vanishes, $dh = 0$, and $\Id = hd$.  Therefore $C$ and $H(C)$ are homotopy equivalent, and $\pi$ and $\iota$ are homotopy inverses.

			To prove the lemma, notice that $\pi \circ f \circ \iota = \Id$ on $H(C)$.  Therefore $f \simeq \iota \circ \pi \simeq \Id$.
		\end{proof}

		To summarize:

		\begin{Cor}\label{cor:conicMeansCondition}
			Suppose $\hkft$ is a conic pretheory which satisfies every condition except \ref{kft:conditions}.  Then $\hkft$ satisfies Condition \ref{kft:conditions}.
		\end{Cor}

	\subsection{Lee-Bar-Natan homology}\label{subsec:LeeTurner}
		Of the strong Khovanov-Floer theories explored in this paper, Lee-Bar-Natan is the closest to Khovanov homology.  See the footnote in the middle of the Introduction for discussion of the name.

		Lee-Bar-Natan homology can be defined in the same style as Khovanov homology.  The only difference is that the merge map $m$ is defined by
		\begin{align*}
			m(v_+ \otimes v_+) &= v_+ \\ 
			m(v_- \otimes v_-) &= m(v_- \otimes v_+) = m(v_+ \otimes v_-) = v_-
		\end{align*}
		and the split map is $\Delta$ defined by
		\begin{align*}
			\Delta(v_+) &= v_+ \otimes v_- + v_- \otimes v_+ + v_+ \otimes v_+ \\
			\Delta(v_-) &= v_- \otimes v_-.
		\end{align*}
		Recall that there is a grading $q$, the \emph{quantum grading}, on $\CBn$.  For a canonical generator $x \in \CBn(\diagram(I))$,
		\[
			q(x) = \tilde{q}(x) + \|I\| + n_+ - 2 n_-
		\]
		where $\|I\|$ is the sum of the coordinates in $I$, $n_+$ and $n_-$ are the numbers of positive and negative crossings in $\diagram$, and $\tilde{q}(x)$ is the sum of the pluses and minuses in $x$.
		\begin{Thm}\label{thm:leeBarNatan}
			Lee-Bar-Natan homology defines a conic, strong Khovanov-Floer theory.
		\end{Thm}
		\begin{proof}
			Define $\hkft_{\mathrm{BN}}(\diagram) = \CBn(\diagram)$, filtered by the quantum grading.  There is no auxiliary data, so Conditions \ref{kft:coherence} and \ref{kft:handles} are satisfied vacuously.  Condition \ref{kft:crossingless} is clear.  (In fact, the condition is satisfied without taking homology.) The proof of Condition \ref{kft:kunneth} is exactly the same as in Khovanov homology.

			Let $\diagram$ be a link diagram.  To a zero-handle attachment $\handle \co \diagram \to \diagram \cup \, \circ$ assign the map 
			\[
				\CBn(\diagram) \to \CBn(\diagram) \otimes R \cong \CBn(\diagram) \otimes \CBn(\circ) \cong \CBn(\diagram \cup \circ)
			\] 
			given by $X \mapsto X \otimes 1$.  Define the linear map $\epsilon \co R \to \bF$ by $\epsilon(X) = 1$ and $\epsilon(1) = 0$.  To a two-handle attachment map $\handle \co \diagram \cup \circ \to \diagram$ assign the map
			\[
				\CBn(\diagram \cup \circ) \cong \CBn(\diagram) \otimes \CBn(\circ) \cong \CBn(\diagram) \otimes R \to \CBn(\diagram)
			\]
			in which the last arrow is given by $x \otimes r \mapsto \epsilon(r)x$.  The zero- and two-handle maps are given by standard unit and counit maps.  Condition \ref{kft:frobenius} is obvious for the Frobenius algebra $\bF[X]/(X^2 + X)$.  Condition \ref{kft:isotopy} is well-known: a planar isotopy between $\diagram$ and $\diagram'$ gives a bijection between closed components of corresponding resolutions of the two diagrams.  This bijection induces a linear map which is obviously a chain isomorphism.
		\end{proof}

		The normal and strong Khovanov-Floer theories for Lee-Bar-Natan homology are quite different.  Consider a twice-punctured surface of genus greater than one embedded in $\br^3$ as a cobordism between two unknots.  It induces a map $\CBn(U) \to \CBn(U)$.  The induced map on homology -- the subject of strong Khovanov-Floer theories -- is an isomorphism, but the associated graded map is zero.
	\subsection{\Szabo's geometric link homology theory}\label{subsec:szabo}

		Like Lee-Bar-Natan homology, \Szabo's geometric link homology theory is defined via cube of resolutions.  \Szabo's theory also includes maps along the higher-dimensional diagonals of the cube.  Each of these diagonals spans some face in the cube.  The diagonals can be represented by drawing, on the diagram at the initial vertex of the face, the traces of all the saddles cobordisms encountered when traveling to the terminal vertex.  \Szabo assigns maps to such pictures.  Below we briefly review this construction from \cite{Szabo2015}, mostly to set vocabulary and notation. 

	 	A crossingless link diagram along with a collection of properly embedded, oriented arcs is called a \emph{configuration}.  A choice of orientations on the arcs is called a \emph{decoration}; we will often conflate the orientation with the oriented arc.  The number of arcs is the \emph{dimension} of the configuration.  

	 	Let $\diagram$ be a link diagram with $k$ crossings and let $I_0$ be the resolution of all zeros.  Draw arcs where the crossings were (going from right to left in Figure \ref{fig:addCrossing}). The orientations on these $k$ arcs induce orientations on every other decoration in the cube.  A cube of resolutions with such decorations is called a \emph{decorated cube of resolutions}.  Write $\Decos(\diagram)$ for the set of decorations of $\diagram$.

	 	\begin{figure}
	 		\centering
	 		\includegraphics[width=.5\linewidth]{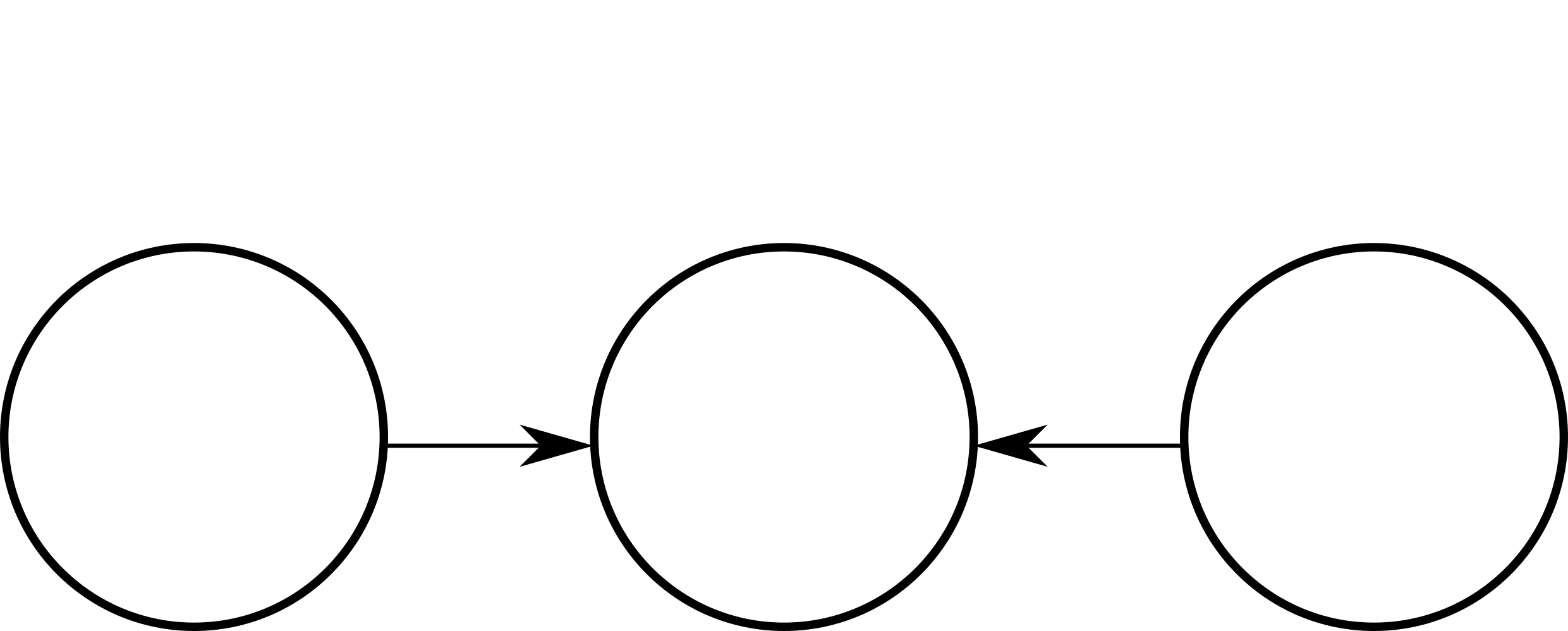}
	 		\caption{A two-dimensional configuration.}
	 		\label{fig:config2}
	 	\end{figure}

	 	Let $\diagram$ be a link diagram with decorations $\decos$.  For any resolution $I$, define 
	 	\[
	 		\CSz(\diagram(I); \decos) = \CKh(\diagram(I)).
	 	\]  
	 	and
	 	\[
	 		\CSz(\diagram; \decos) = \bigoplus_I \CSz(\diagram(I); \decos).
	 	\]  
	 	In other words, the vector space underlying $\CSz$ is the same as that underlying $\CKh$.  Let $I$ and $J$ be resolutions with $I < J$.  If $\|J - I\| = n$, then there is a unique $n$-dimensional face in the cube of resolutions for $\diagram$ with initial vertex $I$ and terminal vertex $J$.  It can be described by an $n$-dimensional configuration $\config(I,J) \co \diagram(I) \to \diagram(J)$.  \Szabo defines a map
	 	\[
	 		\cmap_{\config(I,J); \decos} \co \CSz(\diagram(I); \decos) \to \CSz(\diagram(J); \decos)
	 	\]
	 	based on the topology of the configuration.  The maps assigned to one-dimensional configurations are exactly the maps in Khovanov homology.  Define 
	 	\[
	 		\partial_{i; \decos} = \sum_{\|J - I\| = i} \cmap_{\config(I,J); \decos}
	 	\]
	 	and
	 	\[
	 		\partial_{\decos} = \partial_{1; \decos} + \partial_{2; \decos} + \cdots
	 	\]

	 	\begin{Thm*}[\Szabo, \cite{Szabo2015}]\label{thm:SzabosInvariant}
	 	\begin{enumerate}
	 		\item $(\CSz(\diagram; \decos), \partial_{\decos})$ is a complex.  

	 		\item Let $\decos$ and $\decos'$ be two different sets of decorations for $\diagram$.  Then $(\CSz(\diagram; \decos), \partial_{\decos}) \simeq (\CSz(\diagram; \decos'), \partial_{\decos'})$.  

	 		\item If $\diagram$ and $\diagram'$ are diagrams for the link $L$, then $(\CSz(\diagram; \decos), \partial_{\decos}) \simeq (\CSz(\diagram'; \decos'), \partial_{\decos'})$ for any choices of decorations.
	 	\end{enumerate}
	 	\end{Thm*}

	 	Let $\Sz(\diagram)$ stand for the homology of $\CSz(\diagram; \decos)$.

		\begin{Prop}[Baldwin, Hedden, Lobb \cite{BaldwinHeddenLobb2015}]\label{prop:szaboKunneth}
			$\CSz(\diagram \coprod \diagram') \cong \CSz(\diagram) \otimes \CSz(\diagram')$.  The isomorphism is the same as in the Khovanov K\"{u}nneth formula.
		\end{Prop}
		To prove naturality under choice of decoration, we recall the proof that the chain homotopy type of $\CSz$ is invariant.  Let $\diagram$ be a link diagram with crossings.  Choose one and call it $c$.  Suppose that $I < I'$ are resolutions with differ only at $c$.  Then the configuration $\config(I,I')$ has a single decoration.  Define a linear map
		\[
			H_{I,I'} \co \CSz(\diagram(I)) \to \CSz(\diagram(I')).
		\]
		$H_{I,I'}$ acts as the identity except on the circles affected by the saddle.  If that cobordism merges two circles, then 
		\[
			H_{I,I'}(v_- \otimes v_-) = v_-.
		\]  
		If that cobordism splits two circles, then 
		\[
			H_{I,I'}(v_+) = v_+ \otimes v_+.
		\]  
		Define $H_{I,I'}$ to be zero on any other combination of generators.  Now define $H_c \co \CSz(\diagram; \decos) \to \CSz(\diagram; \decos')$ by
		\[
			H_c = \sum_{\| I' - I \| = 1} H_{I,I'}
		\]
		Define $G_c(x) = \Id + H_c(X)$.  The second part of \Szabo's theorem follows from the observation that $G_c$ is a filtered isomorphism of complexes.  Note that $H_c^2 = 0$ and that the construction of $H_c$ does not depend on the decoration at $c$.  If one were to flip the decoration at $c$ twice, then the isomorphism provided by this construction is
		\[
			(\Id + H_c) \circ (\Id + H_c) = \Id.
		\]
		Observe also that, as a linear map, $H_c$ does not depend on the decoration at any other crossing. 	
		\begin{Lem}\label{lem:naturalityOfSz}
			Let $\diagram$ be a link diagram with at least two crossings.  Call two of them $a$ and $b$.  Choose decorations $\decos$.  Let $\decos_{a}$ be the set of decorations which disagrees with $\decos$ only at $a$.  Define $\decos_b$ and $\decos_{ab}$ similarly.  

			Define
			\begin{align*}
				G_a &\co \CSz(\diagram; \decos) \to \CSz(\diagram; \decos_a) \\
				G_b &\co \CSz(\diagram; \decos) \to \CSz(\diagram; \decos_b) \\
				G_{a,ab} &\co \CSz(\diagram; \decos_a) \to \CSz(\diagram; \decos_{ab}) \\
				G_{b,ab} &\co \CSz(\diagram; \decos_b) \to \CSz(\diagram; \decos_{ab})
			\end{align*}
			Then $G_{b,ab} \circ G_b = G_{a,ab} \circ G_a$.  Therefore 
			\[
				\{\CSz(\diagram; \decos) : \decos \in \Decos(\diagram)\}
			\] 
			and the change of decoration isomorphisms form a transitive system of groups.
		\end{Lem}

		\begin{proof}
			Define $H_a$, $H_b$, $H_{a,ab}$, and $H_{b,ab}$ as above.  We have
			\[
				G_{a,ab} \circ G_a = (\Id + H_b) \circ (\Id + H_a) = \Id + H_a + H_b + (H_b \circ H_a)
			\]
			and
			\[
				G_{b,ab} \circ G_b = (\Id + H_a) \circ (\Id + H_b) = \Id + H_a + H_b + (H_a \circ H_b)
			\]
			So it suffices to show that
			\[
				H_a \circ H_b = H_b \circ H_a.
			\]
			Let $\diagram(I)$ be a decorated resolution in which $a$ and $b$ are both zero-resolved.  (Otherwise $H_a \circ H_b = H_b \circ H_a = 0$.)  Call a  component $a$-\emph{active} (resp. $b$-active) if it intersects the decoration for $a$ (resp. $b$).  If no component is both $a$- and $b$-active, then the $H$ maps clearly commute.  Otherwise, the active components and decorations form one of the two-dimensional configurations in Figure 2 of \cite{Szabo2015}.  There are sixteen such configurations, but as the $H$ maps does not depend on the orientations of the decorations there are effectively only eight.  It is not hard to check that  $H_a \circ H_b = H_b \circ H_a$ in each.
		\end{proof}

		Now we can define diagrammatic handle attachment maps.  Let $(\diagram, \decos)$ be a decorated link diagram.  Observe that if $\diagram'$ is the result of a diagrammatic attachment then $\decos$ defines a set of decorations $\decos'$ on $\diagram'$.  This defines a map
		\[
			\phi \co \Decos(\diagram) \to \Decos(\diagram').
		\]
		Define zero- and two-handle attachment maps as in Khovanov homology.  It follows from the K\"{u}nneth formula that these are chain maps.  Define one-handle attachment maps using the recipe at the top of the section.  $\handle[\gamma]$ is the one-handle map assigned to a \emph{oriented} planar arc $\gamma$.  Concretely, $\handle[\gamma]$ is a sum over all the configurations in the differential of $\CSz(\diagram)$ with $\gamma$ added in, as well as one-dimensional configurations given by $\gamma$.  The map therefore depends on the orientation of $\gamma$.

		Given two oriented, disjoint arcs $\gamma$ and $\gamma'$, consider the composition $\handle[\gamma'] \circ \handle[\gamma]$.  (We will freely abuse notation by conflating $\gamma'$ with its image after handle attachment along $\gamma$ and vice-versa.)  Form the diagram $\diagram_{\gamma'\gamma}$ in which both $\gamma$ and $\gamma'$ are replaced by crossings.  Write $\diagram_{\gamma'\gamma;ij}$ with $i, j \in \{0,1\}$ for the diagrams given by resolving $\gamma'$ and $\gamma$ according to $i$ and $j$, respectively.  Write $g_{\gamma'\gamma}$ for the component of the differential which extends from $\CSz(\diagram_{xy;00}) \to \CSz(\diagram_{xy;11})$ -- this is the sum of all configurations which include both $\gamma$ and $\gamma'$.

		Movie move 15 invariance follows from Corollary \ref{cor:conicMeansCondition}, but it's fun to prove it specifically for \Szabo's theory.

		\begin{Lem}\label{lem:szabo16}
			The handle attachment maps for \Szabo's theory satisfy movie move 15.
		\end{Lem}
		\begin{proof}
			Consider the composition of maps
			\[
				\CSz(\diagram) \to \CSz(\diagram \cup \circ) \to \CSz(\diagram)
			\]
			where the first map is a zero-handle attachment and the second is a saddle connecting the new component with some other part of the diagram.  The zero-handle attachment map sends a canonical generator $x$ to $x \otimes v_+$.  The one-dimensional part of the handle attachment map is exactly the Khovanov map: $x \otimes v_+ \mapsto x$.  Any higher-dimensional configuration which is assigned a non-zero map must be of type E because the new component has degree one, see Section 4 of \cite{Szabo2015}.  Any type E map must vanish on the image of the zero-handle attachment because the degree one component is labeled $v_+$.  Therefore the composition of the two maps is the identity on the nose.

			The other direction of movie move 15 follows from a similar computation (or from the duality rule, Section 2 of \cite{Szabo2015}).
		\end{proof}

	\begin{Thm}\label{thm:Szabo}
		\Szabo's geometric spectral sequence defines a conic strong Khovanov-Floer theory.  It is therefore a functorial link invariant.
	\end{Thm}

	\begin{proof}
		Condition \ref{kft:coherence} is satisfied by Lemma \ref{lem:naturalityOfSz}.  Condition \ref{kft:crossingless} is clear even without taking homology, and Condition \ref{kft:kunneth} is Proposition \ref{prop:szaboKunneth}.  Condition \ref{kft:isotopy} holds by identical argument to Khovanov homology.  That $\CSz$ is conic follows immediately from the definition of the one-handle attachment maps.  Condition \ref{kft:conditions} is the subject of Lemma \ref{lem:szabo16}.

		Now we address Condition \ref{kft:handles}.  Suppose that $\diagram'$ is obtained by a diagrammatic handle attachment on $(\diagram, \decos)$.  We defined a map $\phi \co \Decos(\diagram) \to \Decos(\diagram')$ which is obviously a bijection.  Suppose that $\gamma$ is an oriented arc.  Write $-\gamma$ for the oppositely oriented arc.  We need to show first that $\handle[\decos, \phi(\decos),\gamma] \simeq \handle[\decos, \phi(\decos),-\gamma]$.  Form the diagrams $\diagram_{\gamma}$ and $\diagram_{-\gamma}$.  The complexes $\CSz(\diagram_{\gamma;i})$ and $\CSz(\diagram_{-\gamma;i})$ are identical for $i = 0$ or $1$.  Write $H_\gamma$ for the change of decoration isomorphism at the crossing corresponding to $\gamma$ (and $-\gamma$).  The proof that $\CSz$ is invariant under change of decoration shows that
		\[
			\partial_{\gamma} + \partial_{-\gamma} = H_{\gamma} \circ \partial_{\gamma;0} + \partial_{-\gamma; 1} H_\gamma
		\]
		It follows that $\handle[\decos, \phi(\decos),\gamma] \simeq \handle[\decos, \phi(\decos),-\gamma]$ via the homotopy $H_\gamma$.  The conditions of Lemma \ref{lem:naturality} hold because $H_{\gamma}$ is a chain map on $\diagram_{\gamma}$.

		Condition \ref{kft:handleKunneth} is clear for zero- and two-handles.  For one-handles, the condition holds by the \emph{disconnected rule} which $\cmap$ satisfies.  The \emph{active part} of a configuration $\config$ is the set of connected components of $\config$ which contain a decoration.  A configuration is \emph{disconnected} if its active part is not connected.  The disconnected rule states that if $\config$ is disconnected then $\cmap_{\config} = 0$.  If a diagram splits then any configuration which includes decorations in each part is necessarily disconnected.  Therefore a map $\handle[\gamma]$ with $\gamma$ on only one component must necessarily act as the identity on the other component.
	\end{proof}

\section{Symplectic and gauge-theoretic strong Khovanov-Floer theories}\label{sec:hardExamples}
	\subsection{Heegaard Floer homology of branched double covers}\label{subsec:HF}

		Let $L \subset \br^3$ be a link. Let $\diagram$ be a diagram of $L$ with $k$ crossings.  Let $I^0$ be the all-zeros resolution of $\diagram$.  Draw $k$ small arcs in $\diagram(I^0)$ where the crossings used to be.  Diagrammatic one-handle attachments along these arcs produce all the other resolutions of $\diagram$.

		Write $\Sigma(-L)$ for the double cover of $S^3$ branched along the mirror of $L$.  The small arcs lift to knots in $\Sigma(-\diagram(I_0))$ and handle attachments along them lift to Dehn surgeries.  Form the group
		\[
			X(\diagram) = \bigoplus_{I \in \{0,1\}^k} \CF(\Sigma(-\diagram(I))).
		\]
		(Later will be more specific about what $\CF(-\diagram(I)))$ means.)  Suppose that $I$ and $J$ are resolutions with $I \leq J$.  \Ozsvath and \Szabo \cite{Ozsvath2005b} define a map 
		\[
			d_{I,J} \co \CF(\Sigma(-\diagram(I))) \to \CF(\Sigma(-\diagram(J))).
		\]  
		If $I = J$, then $\partial_{I,J}$ is the usual differential on $\CF(\Sigma(-\diagram(I)))$.  In general, $\partial_{I,J}$ counts certain holomorphic polygons.  Write 
		\[
			\partial = \sum_{I \leq J} \partial_{I,J}.
		\]
		\begin{Thm*}[\Ozsvath, \Szabo \cite{Ozsvath2005b}]
			$(X(\diagram), \partial)$ is a complex.  Its homology is $\HF(\Sigma(-L))$.
		\end{Thm*}

		Observe that $\Sigma(-\diagram(I)) = \#^{m-1} (S^1 \times S^2)$ where $m$ is the number of components of $\diagram(I)$.  It is well-known that $\HF(\#^{m-1}(S^1 \times S^2)) \cong (\bF \oplus \bF)^{\otimes (m-1)}$, so the associated graded object of $(X(\diagram),\partial)$ is isomorphic to $\rCKh(\diagram)$, the reduced Khovanov chain group of $\diagram$.  \Ozsvath and \Szabo show that the length one component of $\partial$ agrees with the Khovanov differential.

		\begin{Thm*}[\Ozsvath, \Szabo \cite{Ozsvath2005b}]
			There is a spectral sequence from $\rKh(\diagram)$ to $\HF(\Sigma(-L))$.  A minor alteration of this construction yields a spectral sequence from $\Kh(\diagram)$ to $\HF(\Sigma(-L) \# S^1 \times S^2) \cong \HF(\Sigma(-L)) \otimes V$.
		\end{Thm*}

		Heegaard Floer homology assigns chain complexes to Heegaard diagrams, not three-manifolds.  At present, there are two consistent systems for building Heegaard diagrams to prove the theorem: the bouquet diagrams used by \Ozsvath and \Szabo, and the branched diagrams developed by the author in \cite{Saltz2017}.  We use the latter because they have a more transparent connection to Khovanov homology and therefore behave better with respect to the Frobenius condition of Definition \ref{def:gkft}
		
		Our goal in this section is to prove the following theorem.

		\begin{Thm}\label{thm:HFisAStrongKFT}
			There is a conic, strong Khovanov-Floer theory called $\hkft_{\HF}$ so that $\hkft_{\HF}(\diagram) = X(\diagram)$.
		\end{Thm}

		Through the rest of this section we will assume familiarity with the fundamentals of Heegaard Floer homology at the level of \cite{Ozsvath2005b}. 

		Let $f \co X(\diagram) \to X(\diagram')$ be a filtered map.  Such a map can be decomposed as $f = f_0 + f_{>0}$ where $f_0$ is the \emph{graded}, i.e. homological degree zero, part of the map.  Call $f_{>0}$ the \emph{higher} part of $f$. Many of the arguments in this section boil down to recognizing that the graded part of some map is well-studied while the higher part vanishes due to the ``Cancellation Lemma'' at the end of the next section.

		\subsubsection{Branched diagrams}\label{subsubsec:branched}
\begin{figure}
				\centering
				\includegraphics[width=.75\linewidth]{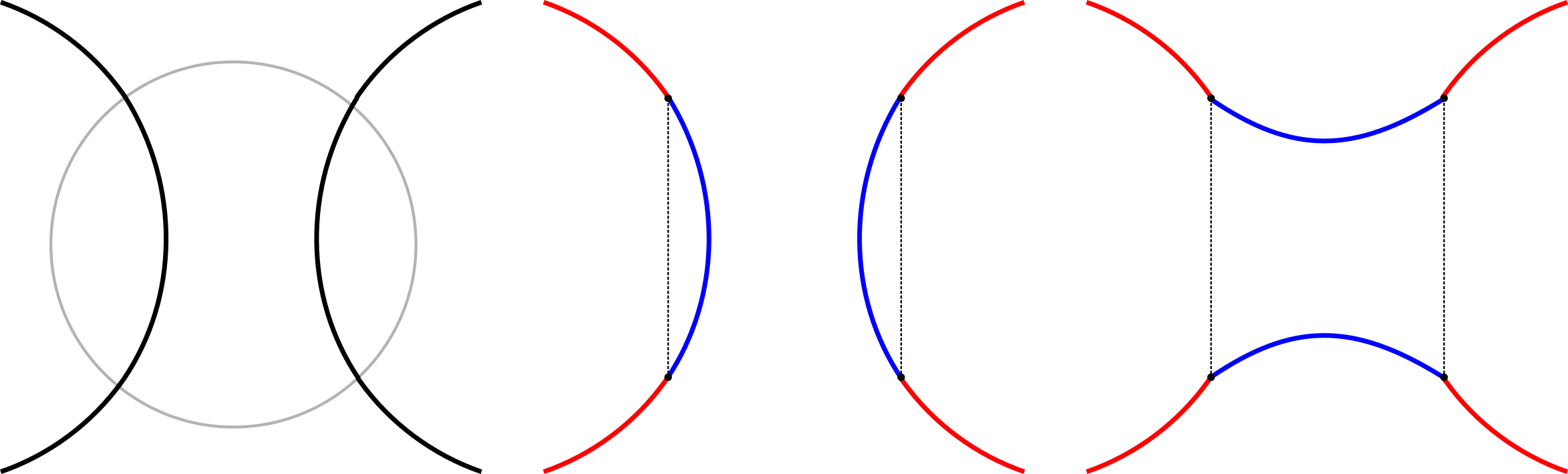}
				\caption{Left: a small circle around what used to be a crossing in $\diagram(I^0)$.  Middle: red and blue curves, as well as small arcs, in $\branched(\diagram(I^0))$.  Right: red and blue curves, as well as small arcs, in $\branched(\diagram(I))$ near a crossing where $I$ and $I^0$ differ.}
				\label{fig:HFbranchedLocal}
			\end{figure}

			Let $\diagram$ be a link diagram with $k$ crossings on a sphere with a basepoint $w \in S^2 \setminus \diagram$.  Draw a little circle around each crossing of $\diagram$ -- this puts $k$ circles on each resolution.  Let $\diagram(I^0)$ be the all-zeros resolution of $\diagram$.  On this diagram, color the arcs outside the circles red and color the arcs inside the circles blue.  Draw small dotted arcs between the red-blue intersection points as shown in Figure \ref{fig:HFbranchedLocal}.  Make a copy of this diagram, calling the second basepoint $w'$.

			If one thinks of the dotted arcs as branch cuts between the two diagrams, then the red arcs form red circles and the blue arcs form blue circles.  Call these circles $\boldA$ and $\boldB(I^0)$, respectively.  These circles lie on a closed surface $V$ with two basepoints $w$ and $w'$.  Write $\mathbf{w} = \{w,w'\}$.  For any other resolution $I$, define $\boldB(I)$ to contain parallel curves to $\boldB(I^0)$ wherever $I$ and $I^0$ agree.  (In this context, ``parallel'' curves are perturbed to intersect twice transversely.) Otherwise, $\boldB(I)$ contains the blue curves on the right of Figure \ref{fig:HFbranchedLocal}.  Define
			\[
				\branched(\diagram(I)) = (V, \boldA, \boldB(I),\mathbf{w}).
			\]
			Call $\branched(\diagram(I))$ the \emph{branched Heegaard diagram} for the branched double cover of $\diagram(I)$.
	
			\begin{Prop*}[\cite{Saltz2017}]
				$\branched(\diagram(I))$ is an admissible Heegaard diagram for the double cover of $S^3$ branched along $\diagram(I)$.  The differential on $\CF(\branched(\diagram(I)))$ vanishes for any choice of analytic data.  There is a unique generator $\Theta_I \in \CF(\branched(\diagram(I))$ with maximal Maslov grading, and it is represented by a single intersection point $\theta_I$. 
			\end{Prop*}

			We will often just write $\CF(\boldA, \boldB(I))$ or just $\CF(I)$.  Similarly, we write $\CF(I,J)$ for $\CF(\boldB(I),\boldB(J))$.  Observe that $\CF(I,J)$ is a diagram for $\#^k(S^1 \times S^2)$ with a unique generator of highest Maslov grading represented by a single intersection point.   

			Let $\mathbf{I} = \{I_1, \ldots, I_n\}$ be a path in the cube of resolutions for $\diagram(I)$; this means that for each $i$, $I_i < I_{i+1}$ and $I_i$ and $I_{i+1}$ differ at only one crossing.  Form the \emph{branched multi-diagram}
			\[
				\branched(\mathbf{I}) = (V, \boldA, \boldB(I_1), \ldots, \boldB(I_n), \mathbf{w}).
			\]
			Let $x \in \CF(I_1)$ and $y_i \in \CF(I_i, I_{i+1})$ for $i \in \{1,n-1\}$ be generators.  Define
			\[
				g_{\mathbf{I}} \co \CF(I_1) \otimes \CF(I_1, I_2) \otimes \CF(I_2, I_3) \otimes \cdots \otimes \CF(I_{n-1},I_n) \to \CF(I_n)
			\]
			by
			\[
				g_{\mathbf{I}}(x \otimes y_1 \otimes \cdots \otimes y_n) = \sum_{z \in \CF(I_n)} \sum_{\substack{\phi \in \pi_2(x, y_1, \ldots, y_n,z) \\ \mu(\phi) = -n}} |\mathcal{M}(\phi)| \, z.
			\]
			Here $\pi_2(x, y_1, \ldots, y_n,z)$ is the set of Whitney $(n+2)$-gons in $\Sym^{g+1}(V)$ connecting those intersection points, $\mu(\phi)$ is the Maslov index of $\phi$, and $|\mathcal{M}(\phi)|$ is the number of holomorphic representatives of $\phi$.  (We use the Maslov index convention from \cite{Manolescu2017}, Section 3.5.)  Define
			\[
				f_{\mathbf{I}} \co \CF(I_1) \to \CF(I_n)
			\]
			by
			\[
				f_{\mathbf{I}}(x) = g_{\mathbf{I}}(x)(x \otimes \Theta_1 \otimes \cdots \otimes \Theta_{n-1})
			\]
			Let 
			\[
				X(\diagram) = \bigoplus_{I \in \{0,1\}^c} \CF(I)
			\]
			and define
			\begin{align*}
				&\partial \co X(\diagram) \to X(\diagram) \\
				&\partial = \sum_{\text{all paths } \mathbf{I}} f_{\mathbf{I}}.
			\end{align*}
			
			\begin{Thm*}[\cite{Saltz2017}]
				$(X(\diagram),\partial)$ is filtered complex.  The associated spectral sequence is isomorphic to \Ozsvath and \Szabo's.  Therefore the second page of the spectral sequence is isomorphic to $\Kh(L)$ and $H(X(\diagram), \partial)$ is isomorphic to $\HF(\Sigma(-L) \# S^1 \times S^2)$.
			\end{Thm*}

			The proof relies on the following fundamental lemma, first proved for bouquet diagrams.  

			\begin{Lem*}[The Cancellation Lemma, \cite{Ozsvath2005b} for bouquet and \cite{Saltz2017} for branched]
				Let $I$ and $J$ be resolutions of $\diagram$.  Let $\mathbf{I}$ be a path from $I$ to $J$.  Write 
				\[
					\tilde{g}_\mathbf{I} \co \CF(I_1, I_2) \otimes \CF(I_2, I_3) \otimes \cdots \otimes \CF(I_{n-1},I_n) \to \CF(I_1,I_n)
				\]
				for the holomorphic $(n+1)$-gon counting map analogous to $g$.  For a path $\mathbf{I}$ from $I$ to $J$ write 
				\[
					\tilde{g}^\Theta_{\mathbf{I}} = \tilde{g}(\theta_1 \otimes \cdots \otimes \theta_{n-1}).
				\]
				Then
				\[
					\sum_{\text{all paths } \mathbf{I}} \tilde{g}^\Theta_{\mathbf{I}} = 0.
				\]
			\end{Lem*}
			
		\subsubsection{Analytic data}

			\begin{Thm}[Baldwin, \cite{Baldwin2011}]\label{thm:HFcomplexInvariance}
				If $X(\diagram)$ and $X(\diagram')$ are complexes built from the same \emph{bouquet} diagrams but with different analytic data, then there is a chain homotopy equivalence $G \co (X,\partial) \simeq (X',\partial')$.
			\end{Thm}

			 The theorem applies to branched diagrams without alteration -- the only appeal to the structure of branched diagrams is to apply the Cancellation Lemma.   Baldwin constructs isomorphisms of spectral sequences which he then bootstraps into filtered chain homotopy equivalences.  This makes it tricky to prove naturality.  The next proof shows how his argument provides chain homotopy equivalences without any reference to spectral sequences.

			\begin{proof}[Proof of Thm \ref{thm:HFcomplexInvariance}]
				Let $\mathbf{I} = \{I_1, \ldots, I_n\}$ be a path in the cube of resolutions for $\diagram$.  For $x \in \CF(I_1)$, define
				\[
					G_{\mathbf{I}}(x) = \sum_{j} g(x \otimes \Theta_1 \otimes \cdots \otimes \Theta_{j-1} \otimes \Theta'_{j-1} \otimes \Theta'_{j} \otimes \cdots \otimes \Theta'_{n})
				\]
				where $\Theta'_k$ means the generator corresponding to $\theta_k$ as an element of $X(\diagram')$ and $g$ counts polygons which are pseudoholomorphic with respect to an interpolating set of analytic data.  Define $G \co X(\diagram) \to X(\diagram')$ by
				\[
					G = \sum_{\text{paths } \mathbf{I}} G_{\mathbf{I}}.
				\]
				Consider the degenerations of the polygons counted by $G$ but with Maslov index one greater.  These `broken polygons' are counted by the sum $G \circ \partial + \partial' \circ G$ (this relies on a cancellation lemma).  They form the boundary of a moduli space of polygons, and therefore their sum must be zero mod two.  So $G$ is a chain map.

				Let $G' \co X(\diagram') \to X(\diagram)$ be the map which interpolates in the opposite direction. A standard degeneration argument shows that 
				\begin{align*}
					G' \circ G - G'_0 \circ G_0 =~ &\partial \circ H + H \circ \partial \\
					&+ \text{(terms which vanish by the cancellation lemma)}
				\end{align*}
				were $H$ counts polygons using analytic data which interpolate from $X(\diagram)$ to $X(\diagram')$ and back to $X(\diagram)$.   Meanwhile, $G'_0 \circ G_0 \simeq \Id$ according to Section 6 of \cite{Ozsvath2004}.  Therefore $G$ and $\widehat{G}$ form a homotopy equivalence with homotopy $H$.
			\end{proof}

			\begin{Prop}\label{prop:HFnaturality}
				Let $X(\diagram)$, $X(\diagram')$, and $X(\diagram'')$ be complexes built from the same branched diagrams but with different analytic data.  Let $G \co X(\diagram) \simeq X(\diagram')$, $G' \co X(\diagram') \simeq X(\diagram'')$, and $G'' \co X(\diagram) \simeq X(\diagram'')$ be the homotopy equivalences of Theorem \ref{thm:HFcomplexInvariance}.  Then $G' \circ G \simeq G''$.
			\end{Prop}

			\begin{proof}
				This time, define $H \co X(\diagram) \to X(\diagram'')$ by
				\[
					H_{\mathbf{I}}(x) = \sum_{j,k} g''(x \otimes \Theta_1 \otimes \cdots \otimes \Theta_j \otimes \Theta'_j  \otimes \cdots \otimes \Theta'_{k} \otimes \Theta''_{k} \otimes \cdots \otimes \Theta''_n )
				\]
				where $g''$ counts polygons which are pseudoholomorphic which respect to a set of analytic data which interpolates between $X(\diagram)$, $X(\diagram')$, and $X(\diagram'')$.

				 Consider the degenerations of the space of polygons counted by $H$ but with Maslov degree one higher.  The degenerations along chords which intersect $\boldA$ correspond to $H \circ \partial$, $\partial'' \circ H$, and $G' \circ G$.  The remaining degenerations cancel by the cancellation lemma (adapted to the situation in which two vertices may be identical, see Baldwin's proof of our Theorem \ref{thm:HFcomplexInvariance} in \cite{Baldwin2011}) except in degenerations of paths in which $j = k$.  In that case, there is a unique degeneration which corresponds to
				 \[
				 	g^*(x \otimes \Theta_1 \otimes \cdots \otimes \Theta_{j-1} \otimes g''(\Theta_j \otimes \Theta'_j \otimes \Theta''_j) \otimes \Theta''_{j+1} \otimes \cdots \otimes \Theta''_n).
				 \]
				 Here $g^*$ interpolates between $X(\diagram)$ and $X(\diagram'')$.  Now $g''(\Theta_j \otimes \Theta'_j \otimes \Theta''_j) = \Theta^*_j$, where $\Theta^* \in \CF(I_j,I''_j)$ as in Theorem 4.5 of \cite{Ozsvath2005b}.   Therefore this degeneration corresponds to $G''$ and $G'' \simeq G' \circ G$.
			\end{proof}

			\begin{Cor}\label{cor:HFisotopy}
				If $\diagram'$ is the result of a planar isotopy of $\diagram$, then $X(\diagram) \simeq X'(\diagram)$.
			\end{Cor}
			\begin{proof}
				A planar isotopy of $\diagram$ induces simultaneous isotopies of all the $\boldA$ and $\boldB$ curves.  This can instead be understood as a simultaneous change in the families of almost-complex structures underlying the construction of $X(\diagram)$.  So Proposition \ref{prop:HFnaturality} implies that the isotopy induces a homotopy equivalence.
			\end{proof}

			For small enough isotopies, the change-of-almost-complex-structure map is the `nearest neighbor' map which sends each generator to its obvious counterpart.

		\subsubsection{Isolated components}\label{subsubsec:isolated}
			The following Proposition will be important in proving that $\hkft$ satisfies Condition \ref{kft:handleKunneth} of Definition \ref{def:gkft}.  

			Let $\mathbf{I} = \{I_1, \ldots, I_n\}$ be a path of resolutions for $\diagram$.  This path corresponds to some sequence of one-handle attachments.  Let $c$ be a component of $\diagram$ which is disjoint from the supports of these one-handles. (Assume that $c$ is not the only component of $\diagram$.)  Call $c$ an \emph{isolated} component. Write $\diagram'(I_j) = \diagram(I_j) \setminus c$.  Then
			\[
				\CF(\diagram(I_j)) \cong \CF(\diagram'(I_j)) \otimes \CF(S^1 \times S^2).
			\]
			\begin{Prop}\label{prop:noOtherPolygons}
				Suppose that $c$ is an isolated component of $\diagram$ with respect to $\mathbf{I}$.  Under the K\"{u}nneth isomorphism above, there is a collection of analytic data so that $f_{\mathbf{I}} \co \CF(\diagram(I_1)) \to \CF(\diagram(I_n))$ acts as the identity on the last factor.
			\end{Prop}
			\begin{proof}
				Let $R$ be a region in $\branched(\diagram'(\mathbf{I}))$ which represents a Whitney polygon of the appropriate Maslov index to be counted by $f'_{\mathbf{I}}$.  Suppose that it connects generators $x$ and $y$.  Write $\Theta_+$ and $\Theta_-$ for the higher and lower Maslov index generators of $\CF(S^1 \times S^2)$.  Figure \ref{fig:HFStretching} shows two regions, $R_+$ and $R_-$, in $\CF(\diagram(\mathbf{I})$ so that $R + R_{\pm}$ is a Whitney polygon from $x \otimes \Theta_\pm$ to $y \otimes \Theta_\pm$.  These regions have exactly the Maslov index so that $R + R_{\pm}$ can be counted by $f_{\mathbf{I}}$.

				At the level of polygons, the region $R$ represents some holomorphic map $\phi \co S \to V'$ where $S$ is a branched cover of the unit disk $\mathbb{D} \subset \mathbb{C}$.  The regions $R_{\pm}$ represent holomorphic maps $\phi_{\pm}$ from $\mathbb{D}$ to $V$.  So $R + R_{\pm}$ represents $\phi \times \phi_{\pm}$, still a maps from a branched cover of $\mathbb{D}$ to $V$.  By the Riemann mapping theorem there are unique holomorphic representatives of $\phi_{\pm}$.  Therefore the operations
				\[
					R \mapsto R + R_{\pm}
				\] 
				induce injective maps 
				\[
					\mathcal{M}(\phi) \to \mathcal{M}(\phi \times \phi_{\pm}).
				\]
				The images of these two maps are clearly distinct.  We claim that the union of their images is all of $\mathcal{M}(\phi \times \phi_{\pm})$.  In other words, for some choice of analytic data, any polygon represented by a region $R'$ which is counted by $f_{\mathbf{I}}$ must be of the form $R + R_+$ or $R + R_-$.  If $R'$ fixes the labeling on $c$, this is the argument of Proposition 7.11(b) in \cite{Manolescu2017}.  The point is to consider families of almost-complex structures parametrized by $s$, $s'$, $t$, and $t' \in \mathbb{R}_{>0}$; the reciprocals of these quantities are the lengths of the curves $\sigma$, $\sigma'$, $\tau$, and $\tau'$, respectively, in Figure \ref{fig:HFStretching}.  (This is ``neck-stretching;'' the curve gets short when its normal neighborhood is stretched in the normal direction.) In the limit as $s, s', t, t' \to \infty$, polygons are stretched into two pieces.  Compactness results imply that there are finite values of $s$ and $t$ so that the moduli space of pseudoholomorphic polygons is diffeomorphic to a space of broken polygons.
				
				\begin{figure}
					\centering
					\def\svgwidth{\linewidth}{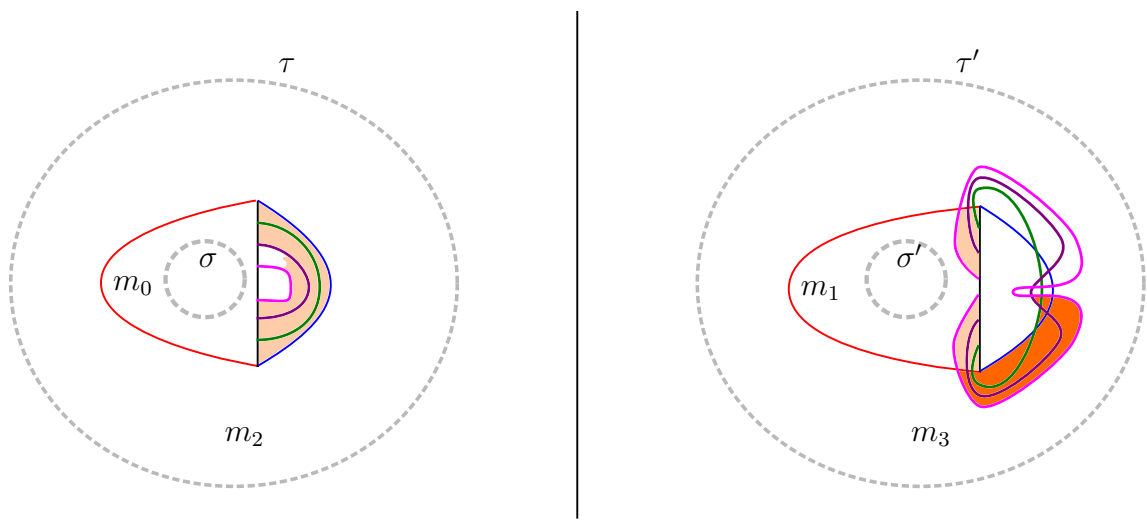}
					\caption{The curves $\tau$, $\tau'$, $\sigma$, and $\sigma'$ around an isolated component.  The dark region is $R_+$.  The dark region together with the light region forms $R_-$.}
					\label{fig:HFStretching}
				\end{figure}

				One difference between our setting and Manolescu-\Ozsvath's is that we are stretching several necks rather than one, but this is not important for the analysis.  The surjectivity result follows.  (Another difference is that they study a connect sum with a spherical Heegaard diagram, while ours is topologically more complicated.  This will only affect the Maslov index computation below.)

				Now we show that every domain counted by $f_\mathbf{I}$ must preserve the labeling on $c$.  Let $R$ be a region which changes the labeling.  Write it as a pair $(R_{\text{small}}, R_{\text{big}})$ where $R_{\text{small}}$ is the part of $R'$ within $\tau$, $\tau'$, $\sigma$, and $\sigma'$ and $R_{\text{big}}$ is the part outside of $\tau$, $\tau'$, $\sigma$, and $\sigma'$.  Identify $R_{\text{big}}$ with the obvious region on $\CF(\diagram'(\mathbf{I}))$.  Observe that, in this region, $m_0$ and $m_2$ count the same multiplicity, and similarly for $m_1$ and $m_3$.  Write $m'_0$ and $m'_1$ for these multiplicities in $R'$.  Then 
				\begin{align*}
					m'_0 &= m_0 + m_2 \\
					m'_1 &= m_1 + m_3.
				\end{align*}
				We can compute $\mu(R)$ in terms of $\mu(R_{\text{small}})$ and $\mu(R_{\text{big}})$:
				\begin{equation}\label{eqn:maslov1}
					\mu(R) = \mu(R_{\text{small}}) + \mu(R_{\text{big}}) - (1-n).
				\end{equation}
				where $n$ is the length of $\mathbf{I}$.  This follows from Sarkar's formula for the Maslov index \cite{Sarkar2011}.   The correction $(1-n)$ is for the last term in Sarkar's formula.  (It may seem odd that there is no correction for the Euler measure, but the comparison of the $m'$s with the $m$s shows that none is necessary.)

				Now we show that
				\begin{equation}\label{eqn:maslov2}
					\mu(R_{\text{small}}) = (1-n) + m_0 + m_1 + m_2 + m_3.
				\end{equation} 
				The formula holds for $R_+$.  Any other domain representing a relevant Whitney polygon can be described by adding and subtracting Whitney disks with Maslov index one on the small side of $V$.  (The resulting domain may not have a pseudoholomorphic representative, but it will have the same Maslov index as any pseudoholomorphic representative in the same class of Whitney polygon.)  The index of such a disk is exactly its contribution to $m_0 + m_1 + m_2 + m_3$, so the formula holds for all small regions.

				Combine equations \ref{eqn:maslov1} and \ref{eqn:maslov2} to obtain
				\[
					\mu(R) = \mu(R_{\text{big}}) + m_0 + m_1 + m_2 + m_3.
				\]
				Now if $R'$ changes the label on $c$, then it must have a corner at the higher $\boldA$-$\boldB(I_1)$ generator.  This implies that $m_0 + m_1 + m_2 + m_3 > 0$.  Therefore
				\[
					\mu(R) > \mu(R_{\text{big}}).
				\]
				If $R_{\text{big}}$ is to be counted in the degenerated setting, then $\mu(R_{\text{big}}) \geq 1 - n$.  Therefore 
				\[
					\mu(R) > 1 - n.
				\]
				But such an $R'$ cannot be counted by $f_{\mathbf{I}}$.
			\end{proof}
	
		\subsubsection{Cobordism maps}
			Suppose that $\diagram'$ is obtained from $\diagram$ by a diagrammatic zero-handle attachment attached far away from $\diagram$.  For every resolution $I$ the diagram $\branched(\diagram'(I))$ looks like $\branched(\diagram(I))$ with the branched diagram of an unknot attached between the basepoint.  So $\CF(\diagram'(I)) = \CF(\diagram(I)) \otimes \CF(S^1 \times S^2)$. The new component is isolated.

			Define the zero-handle attachment map by
			\begin{align*}
				&Z_I \co \CF(\diagram(I)) \to \CF(\diagram'(I)) \\
			 	&Z_I(x) = x \otimes \Theta
			 \end{align*}
			 and
			 \begin{align*}
			 	&Z \co X(\diagram) \to X(\diagram') \\
			 	&Z = \sum_{I} Z_I
			\end{align*}
			
			\begin{Lem}\label{lem:HFZeroHandles}
				$Z$ is a chain map. 
			\end{Lem}
			\begin{proof}
				It is obvious that $Z_I$ is a chain map -- $\CF(I)$ has no differential.  The fact that $Z_I$ commutes with all the polygon counts in $\partial$ is implied by Proposition \ref{prop:noOtherPolygons}.
			\end{proof}

			Define two-handle attachment maps as the duals of the zero-handle attachment maps; i.e. $T \co \CF(\branched(\diagram \cup \circ)) \to \CF(\branched(\diagram))$ is given by 
			\begin{align*}
				T(\mathbf{x} \otimes \Theta) &= 0 \\
				T(\mathbf{x} \otimes \Theta_-) &= \mathbf{x}
			\end{align*}
			where $\CF(\branched(\diagram \cup \circ))$ is identified with $\CF(\branched(\diagram) \otimes \CF(S^1 \times S^2)$.  One-handle attachment maps are defined according to the recipe at the beginning of Section \ref{sec:examples}.

			\begin{proof}[Proof of Theorem \ref{thm:HFisAStrongKFT}]
				Proposition \ref{prop:HFnaturality} proves that Condition \ref{kft:coherence} holds with.  Condition \ref{kft:crossingless} is a standard calculation.  The K\"{u}nneth formula follows from the K\"{u}nneth formula in Heegaard Floer homology, see Section 6 of \cite{Ozsvath2004b}.
				
				To check that the theory satisfies Condition \ref{kft:handles}, one must construct a map which takes analytic data for $X(\diagram)$ to analytic data for $X(\diagram')$.  Connect the two basepoint regions with a tube to obtain singly-pointed Heegaard diagrams $\branched_1$ -- this does not change the represented three-manifold or any of the Heegaard Floer data.  We work with the singly pointed diagrams $\branched_1$ so that the handle attachments correspond to connect summing and its inverse.  Now it is clear that there is a surjective map $\phi$ from data for $X(\diagram')$ to $X(\diagram)$: choose a path from your favorite data $A$ to connect sum data $A' \times A_2$ where $A_2$ is some standard data on a genus two surface.  Then $\phi(A) = A'$.  The space of data is contractible, so there's no obstruction to finding such a path.  The reverse map $A' \mapsto A' \times A_2$ is not surjective, but it still satisfies Lemma \ref{lem:naturality2}.  

				Condition \ref{kft:frobenius} is immediate for zero- and two-handles.  For one-handles, it is straightforward diagrammatic computation.  Alternatively, one can appeal to the equivalence between the branched and bouquet spectral sequences.  Let $\diagram$ be a crossingless link diagram.  Let $\diagram'$ be the one-crossing diagram given by attaching a crossing to $\diagram$ along some planar arc $\gamma$.  Write $X(\diagram')$ for the branched complex and $X'(\diagram')$ for the bouquet complex and $E^i$, $E'^i$ for the corresponding spectral sequences.  There is a map $f \co X(\diagram') \to X'(\diagram')$ which induces a map of spectral sequences $\{f^i\}$.  In \cite{Saltz2017} we showed that $f^i$ is an isomorphism for $i > 0$.  We know that the first page of each spectral sequence is isomorphic to Khovanov homology.  For a one-crossing diagram, there is an isomorphism
				\[
					\CKh(\diagram') \cong E'^1(\diagram') \cong E^1(\diagram')
				\]
				which is equivariant with respect to the differential.  This fixes an isomorphism between $\hkft(\diagram)$ and $\CKh(\diagram)$.  (We leave it to the reader to show that the isomorphism does not depend on the particular added crossing.)  So $d_{\CKh(\diagram')}$ agrees with $\partial_1$ on $X(\diagram)$, and $\partial = \partial_1$ because there's only one crossing.  Therefore the one-handle map agrees with the map on Khovanov homology.  

				Condition \ref{kft:isotopy} follows from Corollary \ref{cor:HFisotopy}.  It suffices to check Condition \ref{kft:handleKunneth} for planar isotopy and each type of handle attachment.  Planar isotopy is handled by $\ref{cor:HFisotopy}$ and the comment following it.  Proposition \ref{prop:noOtherPolygons} handles the handle attachments.  It also shows that maps assigned to handles with disjoint supports commute.
			\end{proof}

	\subsection{Other gauge theories}\label{subsec:gauge}
		Our goal in this section is to suggest proofs that the gauge-theoretic Khovanov-Floer theories are strong.

		In \cite{Kronheimer2011}, Kronheimer and Mrowka assign to a link diagram $\diagram$ with a basepoint $p$ a filtered complex $C(\diagram,p)$.

		\begin{Thm*}[Kronheimer, Mrowka \cite{Kronheimer2011}]
			Suppose that $\diagram$ is a diagram of the link $K$.  Let $p \in K$.  The homology of $C(\diagram,p)$ is the \emph{reduced instanton Floer homology} $I^\natural(K,p)$ of $K$.  The induced spectral sequence has $E_2 = \rKh(K,p)$, the \emph{reduced Khovanov homology of $K$} with basepoint $p$.
		\end{Thm*}

		We will not review the definitions of singular instanton knot homology or reduced Khovanov homology. (The latter is quite brief, see \cite{Khovanov2003}.)  

		It is straightforward to adapt the definition of strong Khovanov-Floer theory to consider a basepoint as long as the underlying Frobenius algebra is $R = \bF[X]/(X^2)$.  (It is possible in other settings, but we defer that to future work.)  A diagrammatic cobordism of pointed diagrams is defined exactly as for unpointed diagrams except that handles must be supported in the complement of a neighborhood of the basepoint.  The results of Section \ref{sec:diagrams} all apply, \emph{mutatis mutandis}, see \cite{BaldwinHeddenLobb2015}.  Given that, the definition is essentially identical to Definition \ref{def:gkft} except for the K\"{u}nneth formula:
	    \[
        	\hkft(\diagram \# \diagram',p'') \simeq \hkft(\diagram,p) \otimes \hkft(\diagram',p')
        \]
        where $p''$ lies on the connected sum region and the connected sum is taken between the arcs containing $p$ to $p'$.  Cobordisms between unlinks are still given by operations in the Frobenius algebra, except that the pointed component acts like $R/(X)$.

		Let $\diagram$ be a pointed link diagram.  Then $I^\natural(\diagram \cup \circ) \cong I^{\natural} \otimes R$ comes with all the naturality one could ask for, see Section 8 of \cite{Kronheimer2011}.  It is shown in Section 2 of \cite{KronheimerMrowka} that this isomorphism can be achieved at chain homotopy level: there is a collection of analytic data so that the chain complex underlying $I^{\sharp}(\circ)$ is isomorphic to $V$ with vanishing differential.  So definitions of zero- and one-handle attachment maps are clear.

		\begin{Thm}\label{thm:instantons}
			There is a conic, pointed homotopy Khovanov-Floer theory $\hkft_I$ which assigns to pointed diagrams $(\diagram,p)$ the complex $C(\diagram,p)$.
		\end{Thm} 

		\begin{proof}		
			Condition \ref{kft:coherence} is discussed in Theorem 3.3 and the beginning of Section 5 of \cite{KronheimerMrowka}.  Condition \ref{kft:crossingless} is Corollary 8.5 of \cite{Kronheimer2011}.  Condition \ref{kft:kunneth} is implied by Corollary 5.9 of \cite{Kronheimer2011}.  

			Condition \ref{kft:handles} is a combination of standard arguments and Theorem 3.3.  Condition \ref{kft:frobenius} follows from Lemma 8.7 of \cite{Kronheimer2011} as well as the definition of the zero- and one-handle maps.  Condition \ref{kft:isotopy} is Proposition 8.10 of \cite{Kronheimer2011}.  Condition \ref{kft:handleKunneth} is Proposition 5.8 of \cite{Kronheimer2011}.  Now Condition \ref{kft:conditions} follows from Corollary \ref{cor:conicMeansCondition}.
		\end{proof}

		The proof for monopole Floer homology should be essentially the same.  The construction of a conic theory is due to Bloom \cite{Bloom2011}.  Naturality should follow by applying an argument analogous to that of our Proposition \ref{prop:HFnaturality} to Theorem 7.4 of \cite{Bloom2011}.

\bibliographystyle{plain}
\bibliography{bibliography.bib}

\end{document}